\documentclass[times,sort&compress,3p]{elsarticle}
\journal{Journal of Multivariate Analysis}
\usepackage[labelfont=bf]{caption}

\usepackage{amsmath,amsfonts,amssymb,amsthm,booktabs,color,csquotes,epsfig,graphicx,hyperref,url}

\theoremstyle{plain}% Theorem-like structures provided by amsthm.sty
\newtheorem{theorem}{Theorem}

\newtheorem{lemma}{Lemma}
\newtheorem{corollary}{Corollary}

\theoremstyle{definition}

\newtheorem{remark}{Remark}
\newtheorem{example}{Example}

\allowdisplaybreaks

\usepackage{dsfont}         % for indicator function \mathds{1}
\usepackage[T1]{fontenc}

\newcounter{aaa}[theorem]

\newcounter{one}[theorem]

\newcommand{\I}{{\mathbb I}}

\newcommand{\R}{{\mathbb R}}

\newcommand{\ppp}{{\bf p}}
\newcommand{\qqq}{{\bf q}}
\newcommand{\uuu}{{\bf u}}
\newcommand{\vvv}{{\bf v}}

\newcommand{\xxx}{{\bf x}}
\newcommand{\UUU}{{\bf U}}
\newcommand{\XXX}{{\bf X}}
\begin{document}

\begin{frontmatter}

\title{Quantifying directed dependence via dimension reduction}
%\title{A bivariate copula capturing the dependence of a random variable \\ and a random vector, its estimation and applications}

\author[1]{Sebastian Fuchs\corref{mycorrespondingauthor}}

\address[1]{Universit{\"a}t Salzburg, Salzburg, Austria}

\cortext[mycorrespondingauthor]{Corresponding author. Email address: \url{sebastian.fuchs@plus.ac.at}}

%Um diesem Problem zu entgehen, soll eine Dimensionsreduktion durchgeführt werden, welche die grundlegende Information über die gerichtete Abhängigkeit der involvierten Größen unberührt lässt. 

\begin{abstract}
Studying the multivariate extension of copula correlation yields a dimension reduction principle, which turns out to be strongly related with the `simple measure of conditional dependence' $T$ recently introduced by \citet{chatterjee2021}.
In the present paper, we identify and investigate the dependence structure underlying this dimension reduction principle,
provide a strongly consistent estimator for it, and demonstrate its broad applicability. 
For that purpose, we define a bivariate copula capturing the scale-invariant extent of dependence of an endogenous random variable $Y$ on a set of $d \geq 1$ exogenous random variables $\XXX = (X_1, \dots, X_d)$,
and containing the information whether $Y$ is completely dependent on $\XXX$, and
whether $Y$ and $\XXX$ are independent.
The dimension reduction principle becomes apparent insofar as the introduced bivariate copula can be viewed as the distribution function of two random variables $Y$ and $Y^\prime$ sharing the same conditional distribution and being conditionally independent given $\XXX$. 
Evaluating this copula uniformly along the diagonal, i.e., calculating Spearman's footrule, 
leads to an unconditional version of Azadkia and Chatterjee's `simple measure of conditional dependence' $T$. 
On the other hand, evaluating this copula uniformly over the unit square, i.e., calculating Spearman's rho, 
leads to a distribution-free coefficient of determination (also known as `copula correlation').
Several real data examples illustrate the importance of the introduced methodology.
\end{abstract}

\begin{keyword} %alphabetical order
Copulas \sep 
Copula correlation \sep
Directed dependence \sep
Markov product \sep
Spearman's footrule \sep
Spearman's rank correlation.

\MSC[2020] Primary 62H20 \sep
Secondary 62E10
\end{keyword}

\end{frontmatter}

%%%%%%%%%%%%%%%%%%%%%%%%%%%%%%%%%%%%%%%%%%%%%%%%%%%%%%%%%%%%%%%%%%%%%%%%%%%%%%%%%%%%%%%%%%%%%%%%%%%%%%%%%%%%%%%%%%%
%%%%%%%%%%%%%%%%%%%%%%%%%%%%%%%%%%%%%%%%%%%%%%%%%%%%%%%%%%%%%%%%%%%%%%%%%%%%%%%%%%%%%%%%%%%%%%%%%%%%%%%%%%%%%%%%%%%
%%%%%%%%%%%%%%%%%%%%%%%%%%%%%%%%%%%%%%%%%%%%%%%%%%%%%%%%%%%%%%%%%%%%%%%%%%%%%%%%%%%%%%%%%%%%%%%%%%%%%%%%%%%%%%%%%%%
\section{Introduction}

Detecting statistical association among several random variables is an ubiquitous task. 
Measures of association capture the many facets of dependence relationships - from classical linear correlation coefficients to indices for detecting monotone associations, tail dependence, asymmetric dependence etc. 
Most concepts of statistical association, however, consider association as an undirected property, 
i.e., the association among the variables remains unchanged when permuting the variables.
This includes Pearson's correlation coefficient, 
measures of concordance (see, e.g., \cite{ccr1996, duf2006, sfx2016moc, gene2011, joe1990, schmidt2006, tay2016}), 
distance covariance and distance multivariance (see, e.g., \cite{bottcher2019, szekely2007}),
the maximal information coefficient (see, e.g., \cite{reshef2011}),
the randomized dependence coefficient (see, e.g., \cite{lopez2013}),
coefficients of tail dependence (see, e.g., \cite{joe1997}), and
various measures of asymmetry (see, e.g., \cite{stadtmuller2014, kamnitui2018, nel2007}).

In manifold situations one variable may have a stronger influence on another variable than vice versa.
This is why quantifying the degree of predictability or explainability of an endogenous random variable $Y$ 
using the information contained in a set of $d \geq 1$ exogenous random variables $\XXX = (X_1, \dots, X_d)$ requires measures of directed dependence.
An index measuring such a degree should be capable of detecting perfect dependence also known as complete dependence (see \cite{fgwt2021, lancaster1963, wt2011}):
$Y$ is said to be completely dependent on $\XXX$ if there exists a measurable function $f$ such that $Y = f(\XXX)$ almost surely.

Quite recently, \citet{chatterjee2021} introduced their so-called `simple measure of conditional dependence' $T$ given (in its unconditional form) by
\begin{equation} \label{Chatterjee.T.Intro}
  T(Y,\XXX)
	:= \frac{\int_{\mathbb{R}} {\rm var} (P(Y \geq y \, | \, \XXX)) \; \mathrm{d} P^{Y}(y)}
					{\int_{\mathbb{R}} {\rm var} (\mathds{1}_{\{Y \geq y\}}) \; \mathrm{d} P^{Y}(y)},
\end{equation}
which is based on \cite{chatterjee2020, siburg2013}, 
quantifies the scale-invariant extent of dependence of $Y$ on $\XXX$, 
and attracted a lot of attention in the past two years; see, e.g., \cite{auddy2021,bernoulli2021,bickel2020,bickel2022,deb2020,deb2020b,shi2021normal,strothmann2022}. 
As $T$ equals $1$ if and only if $Y$ is completely dependent on $\XXX$,
and $T$ is $0$ if and only if $Y$ and $\XXX$ are independent,
$T$ belongs to a class of indices capable of detecting complete dependence and independence which also include \cite{fswt2015, fgwt2021, fgwt2020, wt2011}.
In order to estimate $T$, Azadkia and Chatterjee developed a nearest neighbour based estimation procedure that allows for a dimension reduction.

In the present paper we identify and investigate the bivariate dependence structure underlying this dimension-reducing estimation principle.
We therefore extend the $(d+1)$-dimensional random vector $(\XXX,Y)$ with continuous marginal distribution functions and connecting copula $A$ and consider the $(d+2)$-dimensional random vector $(\XXX,Y,Y^\prime)$ 
with $Y$ and $Y^\prime$ sharing the same conditional distribution and being conditionally independent given $\XXX$, so that the joint distribution function of $(Y,Y^\prime)$ fulfills
$$
  P (Y \leq y, Y^\prime \leq y^\prime)
	= \int_{\mathbb{R}^d} P (Y \leq y \, | \, \XXX=\xxx) \, P(Y^\prime \leq y^\prime \, | \, \XXX=\xxx) \; \mathrm{d} P^\XXX (\xxx).
$$
The copula of $(Y,Y^\prime)$ we then identify as the outcome of the dimension-reducing transformation $\psi$ given by  
$$
  \psi(A)(s,t) := P (G(Y) \leq s, G(Y^\prime) \leq t)
$$
where $G$ denotes the distribution function of $Y$.
The dimension reduction principle becomes apparent insofar as $\psi$ transforms the $(d+1)$-dimensional copula $A$ (corresponding to $(\XXX,Y)$) to the bivariate copula $\psi(A)$ (corresponding to $(Y,Y^\prime)$) but preserves the key information about the directed dependence of the variables $Y$ and $\XXX$,
i.e., the bivariate copula $\psi(A)$ itself contains the information
(i) whether $Y$ is completely dependent on $\XXX$, and
(ii) whether $Y$ and $\XXX$ are independent.

It is straightforward to show that the index $T$ possesses a representation in terms of the well-known dependence measure Spearman's footrule (see, e.g., \cite{sfx2019spearman,genebg2010,nelsen2006})
$$
  T(Y,\XXX)
	%= \phi(\psi(A))
	= 6 \, \int_{0}^1 P \big(G(Y) \leq t, G(Y^\prime) \leq t\big) \; \mathrm{d}t - 2
	= 6 \, \int_{0}^1 \psi(A)(t,t) \; \mathrm{d}t - 2.
$$
Thus, the initial $(d+1)$-dimensional problem of quantifying the scale-invariant extent of dependence of $Y$ on $\XXX$ reduces to a $2$-dimensional one.
Considering $T$ as a map of $\psi(A)$ suggests evaluating $\psi(A)$ also via other measures of bivariate dependence, 
leading to new indices with which independence, complete dependence and also other facets of dependence between $Y$ and $\XXX$ can be quantified:
For instance, evaluating $\psi(A)$ via
\begin{enumerate}[1.]
\item Spearman's rank correlation (Spearman's rho) leads to a distribution-free coefficient of determination $R^2$ which provides a benchmark for the proportion of variance that can be explained in a copula-based (regression) model.
\item Gini's gamma leads to a measure of `indifference' $Q$ being able to detect whether the dependence structures of $(\XXX,Y)$ and $(\XXX,-Y)$ coincide, a kind of `lack of association' property.
\end{enumerate}{}
For either of these resulting dependence measures, we investigate their ability to identify certain facets of dependence between $Y$ and $\XXX$,
show their invariance with respect to a variety of transformations of the exogenous variables $X_1, \dots, X_d$,
and verify the so-called information gain (in-)equality.
Even though we restrict our investigation to the measures of concordance Spearman's footrule, Spearman's rho and Gini's gamma, 
all leading to meaningful measures of dependence, $\psi(A)$ can in principle be applied to any bivariate measure of concordance / dependence.

For the bivariate copula $\psi(A)$ we propose an estimator whose form is reminiscent of the empirical copula, 
but which is actually based on the graph-based estimator of $T$ developed by \citet{chatterjee2021}.
By applying the tools given in \cite{chatterjee2021} we show that the copula estimator is strongly consistent 
from which strong consistency of the plug-in estimators of Spearman's footrule (coinciding with $T$), Spearman's rho (coinciding with the distribution-free $R^2$) and Gini's gamma (coinciding with the measure of `indifference' $Q$) can be derived.

The rest of this contribution is organized as follows: 
In Section \ref{Sect.phi}, we formally define the dimension-reducing transformation $\psi$ mapping every $(d+1)$-dimensional copula $A$ to an exchangeable bivariate one, 
and show that the resulting bivariate copula $\psi(A)$ captures independence and complete dependence. 
%is invariant with respect to a variety of transformations of the original dependence structure $A$ and satisfies the information gain inequality. 
In Section \ref{Sec.Appl} we then apply the resulting copula $\psi(A)$ to (the measures of bivariate dependence) Spearman's footrule, Spearman's rho and Gini's gamma and demonstrate their ability to identify certain facets of dependence.
Section \ref{Sect.Estimation} is devoted to the estimation of $\psi(A)$. 
We prove strong consistency of the proposed estimator and illustrate its small/moderate sample performance. 
Finally, the potential and importance of the dimension-reducing methodology is illustrated by analyzing several real data examples regarding feature selection, copula-based regression analysis and the detection of reflection invariant dependence structures (Section \ref{Sec.DataEx}).
The proofs can be found in the Appendix.

\bigskip
Throughout this paper we will write $\I := [0,1]$ and let $d\geq 1$ denote an integer which will be kept fixed.
Bold symbols will be used to denote vectors, e.g., $\mathbf{x}=(x_1,\ldots,x_d) \in \mathbb{R}^d$.  
The $d$-dimensional Lebesgue measure will be denoted by $\lambda^d$, in case of $d=1$ we simply write $\lambda$.
We will let $\mathcal{C}^{d+1}$ denote the family of all $(d+1)$-dimensional copulas, 
$M$ will denote the comonotonicity copula, 
$\Pi$ the independence copula
and, for $d=1$, $W$ will denote the countermonotonicity copula
(we omit the index indicating the dimension since no confusion will arise).  
For every $A \in \mathcal{C}^{d+1}$ the corresponding probability measure will be denoted by $\mu_A$, 
i.e., $\mu_A([{\bf 0},\uuu] \times [0,v]) = A(\uuu,v)$ for all $(\uuu,v) \in \I^d \times \I$;
for more background on copulas and copula measures we refer to \cite{fdsempi2016,nelsen2006}. 
For every metric space $(\Omega,\delta)$ the Borel $\sigma$-field on $\Omega$ will be denoted by $\mathcal{B}(\Omega)$.

For a copula $A \in \mathcal{C}^{d+1}$ we denote by $A^L$ the $|L|$-dimensional marginal of $A$ with respect to the coordinates in $L \subseteq \{1,\dots,d+1\}$, 
and for the $l$-dimensional marginal of $A$ ($l \geq 2$) with respect to the first $l \in \{1,\dots,d\}$ coordinates we simply write
$A^{1:l} := A^{\{1,\dots,l\}}$.

According to \cite[Theorem 3.4.3]{fdsempi2016} and due to disintegration, 
every copula $A$ fulfills
$$
  A(\uuu,v)
	= \int_{[{\bf 0},\uuu]} K_A(\ppp,[0,v]) \; \mathrm{d} \mu_{A^{1:d}}(\ppp),
$$
where $K_A$ is (a version of) the Markov kernel of $A$:
A Markov kernel from $\mathbb{I}^d$ to $\mathcal{B}(\mathbb{I})$ is a mapping 
$K: \mathbb{I}^d\times\mathcal{B}(\mathbb{I}) \rightarrow \mathbb{I}$ such that for every fixed 
$F\in\mathcal{B}(\mathbb{I})$ the mapping 
$\uuu\mapsto K(\uuu,F)$ is measurable and for every fixed $\uuu\in\mathbb{I}^d$ the mapping 
$F\mapsto K(\uuu,F)$ is a probability measure. 
Given a random vector $\UUU$ with uniformly distributed univariate marginals and a uniformly distributed random variable $V$ on a probability space $(\Omega, \mathcal{A}, P)$ 
we say that a Markov kernel $K$ is a regular conditional distribution of $V$ given $\UUU$ if 
$K (\UUU(\omega), F) = P( V \in F \,|\, \UUU ) (\omega) $ holds $P$-almost surely for every $F\in \mathcal{B}(\mathbb{I})$. 
It is well-known that for each such random vector $(\UUU,V)$ a regular conditional distribution 
$K(.,.)$ of $V$ given $\UUU$ always exists 
and is unique for $P^\UUU$-a.e. $\uuu\in\mathbb{I}^d$,
where $P^\UUU$ denotes the push-forward of $P$ under $\UUU$.
For more background on conditional expectation and general disintegration we refer to \cite{kallenberg2002, klenke2008};
for more information on Markov kernels in the context of copulas we refer to 
\cite{fdsempi2016, sfx2021weak, sfx2021vine}.

We say that $Y$ is completely dependent on $\XXX$ if there exists a measurable function $f$ such that $P(Y = f(\XXX)) = 1$.
If the marginal distribution functions $F_i$ of $X_i$, $i \in \{1,\dots,d\}$, and $G$ of $Y$ are continuous with connecting copula $A$, then (denoting $U_i:=F_i(X_i)$, $i \in \{1,\dots,d\}$ and $V:=G(Y)$) the following statements are equivalent (see \cite{fgwt2021}):
\begin{enumerate}[(i)]
\item
$Y$ is completely dependent on $\XXX$;
\item
There exists a $\mu_{A^{1:d}}-\lambda$ preserving transformation $h: \I^d \to \I$ (i.e., $(\mu_{A^{1:d}})^h = \lambda$)
such that $V = h(\UUU)$ a.s.;
\item
There exists a $\mu_{A^{1:d}}-\lambda$ preserving transformation $h: \I^d \to \I$ 
such that
$ K(\uuu,F)
	:= \mathds{1}_F (h(\uuu)) $
is a regular conditional distribution of $A$. 
\end{enumerate}
For more properties of complete dependence we refer to \cite{lancaster1963} as well as to \cite{fswt2015} and the references therein.

%%%%%%%%%%%%%%%%%%%%%%%%%%%%%%%%%%%%%%%%%%%%%%%%%%%%%%%%%%%%%%%%%%%%%%%%%%%%%%%%%%%%%%%%%%%%%%%%%%%%%%%%%%%%%
%%%%%%%%%%%%%%%%%%%%%%%%%%%%%%%%%%%%%%%%%%%%%%%%%%%%%%%%%%%%%%%%%%%%%%%%%%%%%%%%%%%%%%%%%%%%%%%%%%%%%%%%%%%%%
%%%%%%%%%%%%%%%%%%%%%%%%%%%%%%%%%%%%%%%%%%%%%%%%%%%%%%%%%%%%%%%%%%%%%%%%%%%%%%%%%%%%%%%%%%%%%%%%%%%%%%%%%%%%%
\section{One bivariate copula to capture it all}
%\section{A bivariate copula capturing independence and complete dependence between a random variable and a random vector}
\label{Sect.phi}

Consider a $(d+1)$-dimensional random vector $(\XXX,Y)$ with continuous marginal distribution functions $F_i$ of $X_i$, $i \in \{1,\dots,d\}$, and $G$ of $Y$ and connecting copula $A$.
Then (denoting $U_i:=F_i(X_i)$, $i \in \{1,\dots,d\}$, and $V:=G(Y)$) we have
$(\UUU,V) \sim A$. 
In what follows, we construct an exchangeable bivariate copula capturing the scale-invariant extent of dependence of $Y$ on the random vector $\XXX$ in the sense that it allows to detect (see Theorem \ref{phi.CompleteDep} below)
\begin{enumerate}
\item[(i)]
whether $Y$ and $\XXX$ (or, equivalently, $V$ and $\UUU$) are independent;
\item[(ii)]
whether $Y$ is completely dependent on $\XXX$ (or, equivalently, $V$ is completely dependent on $\UUU$).
\end{enumerate}

\begin{figure}[h!]
		\centering
		\includegraphics[width=0.65\textwidth]{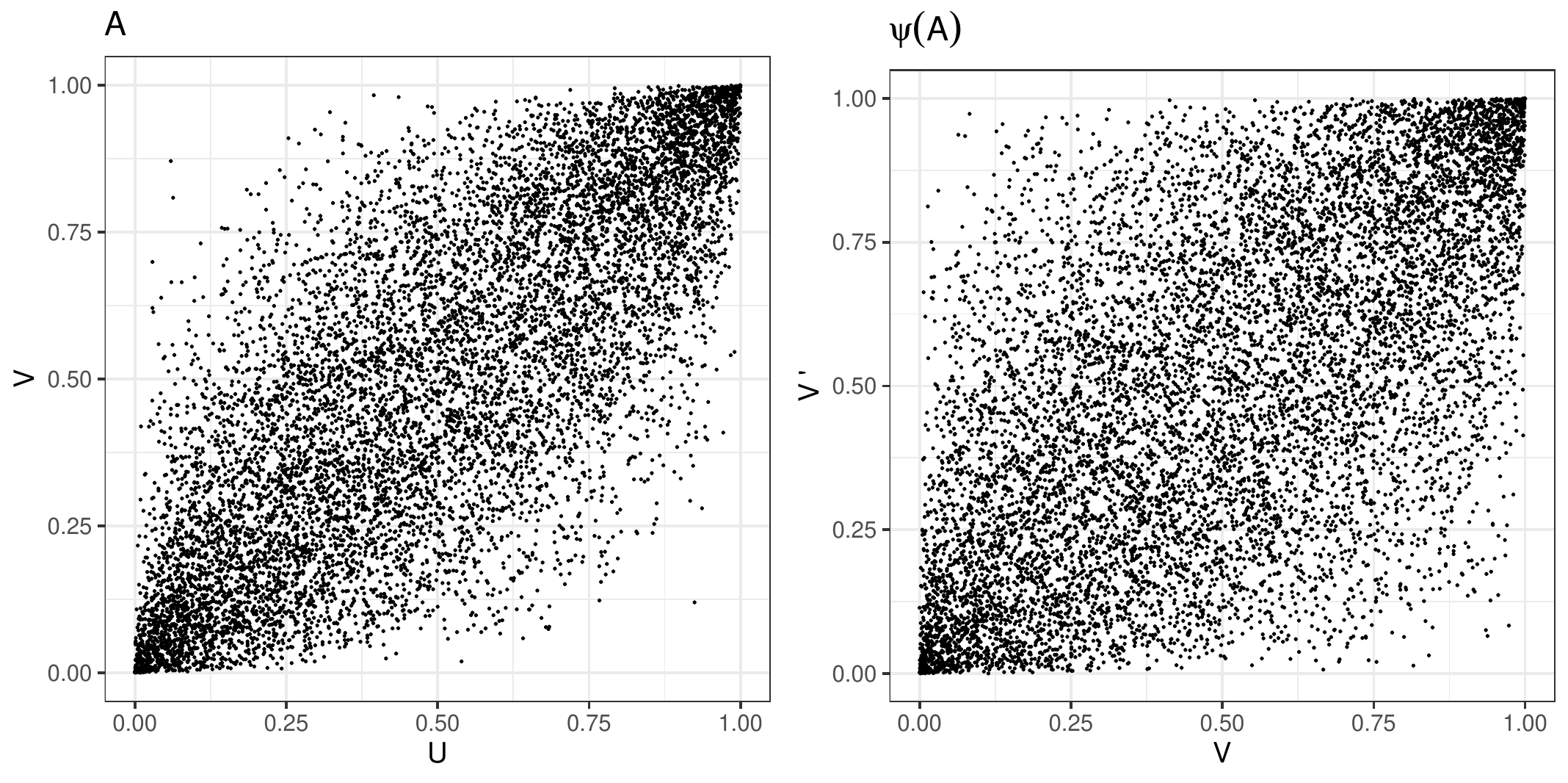}
		\caption{Sample of size $n=10,000$ drawn from a bivariate ($d=1$) Gaussian copula $A$ with correlation parameter $r=0.8$ (left panel), 
together with a sample of equal size from $\psi(A)$ (right panel); due to Example \ref{phi.GaussianII}, $\psi(A)$ is Gaussian with correlation parameter $r^\ast (1) = r^2 = 0.64$.}
		\label{Fig.GaussianII}
\end{figure}

To this end, we extend the random vector $(\UUU,V)$ and consider the $(d+2)$-dimensional random vector $(\UUU,V,V^{\prime})$ with 
$V$ and $V^{\prime}$ sharing the same conditional distribution and being conditionally independent given $\UUU$.
Then, using disintegration, the distribution function of $(V,V^{\prime})$ is a copula and can be expressed as
\begin{eqnarray} 
  P (V \leq s, V^{\prime} \leq t)
	& = &   E \big( P (V \leq s, V^{\prime} \leq t \, | \, \UUU) \big) \label{Phi(A).Representation}
	\; = \; E \big( P (V \leq s \, | \, \UUU) \, P (V^{\prime} \leq t \, | \, \UUU) \big) 
	\\
	& = & \int_{\I^d} P (V \leq s \, | \, \UUU=\uuu) \, P (V^{\prime} \leq t \, | \, \UUU=\uuu) \; \mathrm{d} \mu_{A^{1:d}} (\uuu) \notag
	\\
	& = & \int_{\I^d} K_A (\uuu, [0,s]) \, K_A (\uuu, [0,t]) \; \mathrm{d} \mu_{A^{1:d}} (\uuu), \notag 
	\qquad (s,t) \in \I^2,
	%\\
	%& = & \Pi(s,t) + {\rm cov} (K_A (\UUU, [0,s]), K_A (\UUU, [0,t]))
	%\\*
	%& = & \psi(A) (s,t)
\end{eqnarray}
where $\mu_{A^{1:d}}$ denotes the copula measure of the $d$-dimensional marginal copula $A^{1:d}$ and 
$K_A: \I^d \times \mathcal{B}(\I) \rightarrow \I$ denotes the Markov kernel of $A \in \mathcal{C}^{d+1}$ (with respect to the first $d$ coordinates).
The map $\psi: \mathcal{C}^{d+1} \to \mathcal{C}^2$ given by 
$$ 
  \psi(A) (s,t)
	:= P (V \leq s, V^{\prime} \leq t)
	 = \int_{\I^d} K_A (\uuu, [0,s]) \, K_A (\uuu, [0,t]) \; \mathrm{d} \mu_{A^{1:d}} (\uuu),	 
	\qquad (s,t) \in \I^2,
$$
then transforms every $(d+1)$-dimensional copula to an exchangeable bivariate copula
(see Fig. \ref{Fig.GaussianII} and \ref{Fig.MO} for an illustration), 
but preserves the key information about the directed dependence of the variables involved.

\begin{remark}{}
The copula $\psi(A)$ may be interpreted as a generalization of the well-known bivariate Markov product of copulas: 
For two copulas $B,C \in \mathcal{C}^2$ the \emph{Markov product} $B \ast C: \I^2 \to \I$ given by
$$
  (B \ast C) (s,t)  
	:= \int_{\I} K_{B^\top} (u,[0,s]) \, K_C(u,[0,t]) \; \mathrm{d} \lambda(u)
$$
is a copula (see \cite{darsow1992,wt2012} and \cite[Chapter 5]{fdsempi2016}); here, $B^\top$ denotes the transpose of $B$.
In the bivariate case, i.e., for $d=1$, we therefore have $\psi(A) = A^\top \ast A$. 
Notice that, in this case, the fixed points of $\psi$ are exactly those copulas that are idempotent with respect to the Markov product (i.e., those copulas $A$ satisfying $A \ast A = A$).
\end{remark}

In the remainder of this section, 
we present several key properties of the dimension-reducing transformation $\psi$, 
and we provide explicit formulas for $\psi(A)$ in the case $A$ belongs to certain copula families.

The following theorem shows that the copula $\psi(A)$ characterizes complete dependence of $Y$ on $\XXX$
and independence of $Y$ and $\XXX$; notice that Theorem \ref{phi.CompleteDep} generalizes Theorem 11.1 in \citet{darsow1992}.
For a copula $C \in \mathcal{C}^2$, we denote by $\delta_C$ its diagonal, i.e., $\delta_C(t) := C(t,t)$ for all $t \in \I$.

\begin{theorem}{} \label{phi.CompleteDep}\label{phi.InDep}
Consider a $(d+1)$-dimensional random vector $(\XXX,Y)$ with continuous marginals and connecting copula $A$.
\begin{enumerate}[(i)]
\item
The following statements are equivalent:
\begin{enumerate}[(a)]
\item
$Y$ is completely dependent on $\XXX$;

\item
$\psi(A) = M$.
\end{enumerate}

\item
The following statements are equivalent:
\begin{enumerate}[(a)]
\item
$Y$ and $\XXX$ are independent;

\item
$\psi(A) = \Pi$;

\item
$\delta_{\psi(A)} = \delta_{\Pi}$.
\end{enumerate}
\end{enumerate}
\end{theorem}{}

We now apply the transformation $\psi$ to some well-known parametric copula families:
the class of equicorrelated Gaussian copulas, the class of Marshall-Olkin extreme-value copulas, the class of Fr{\'e}chet copulas
and the class of EFGM copula.
It turns out that $\psi$ transforms every $(d+1)$-dimensional equicorrelated Gaussian copula to a bivariate Gaussian copula where the correlation parameter now depends on the dimension of the conditioning random vector $\XXX$.
We further show that also the class of Marshall-Olkin copulas, 
the class of Fr{\'e}chet copulas and the class of EFGM copulas are to some extent closed with respect to $\psi$.

\begin{example}{} \label{phi.GaussianII}\label{phi.MO.Ex}   \leavevmode
\begin{enumerate}[1.]
\item (Equicorrelated Gaussian copula) 
Consider the $(d+1)$-dimensional equicorrelated Gaussian copula $A_r$ with correlation paramenter $r \in (-1/d,1)$.
After straightforward but tedious calculation, one then finds that $\psi(A_r)$ is a (bivariate) Gaussian copula with correlation paramenter
$$
  r^\ast (d) = \frac{d r^2}{1+(d-1)r}.
$$
Notice that $r^\ast(d) = r \, \frac{d r}{1+(d-1)r} < r$ and, for every $r \in [0,1)$, $\lim_{d \to \infty} r^\ast(d) = r$.
Fig. \ref{Fig.GaussianII} depicts a sample of size $n=10,000$ drawn from a bivariate ($d=1$) Gaussian copula $A$ with correlation parameter $r=0.8$ together with a sample of equal size from $\psi(A)$,
and Fig. \ref{Fig.Gaussian.Cor} illustrates the behaviour of $r^\ast (d)$ when the dimension $d$ increases.
\begin{figure}[h!]
		\centering
		\includegraphics[width=0.65\textwidth]{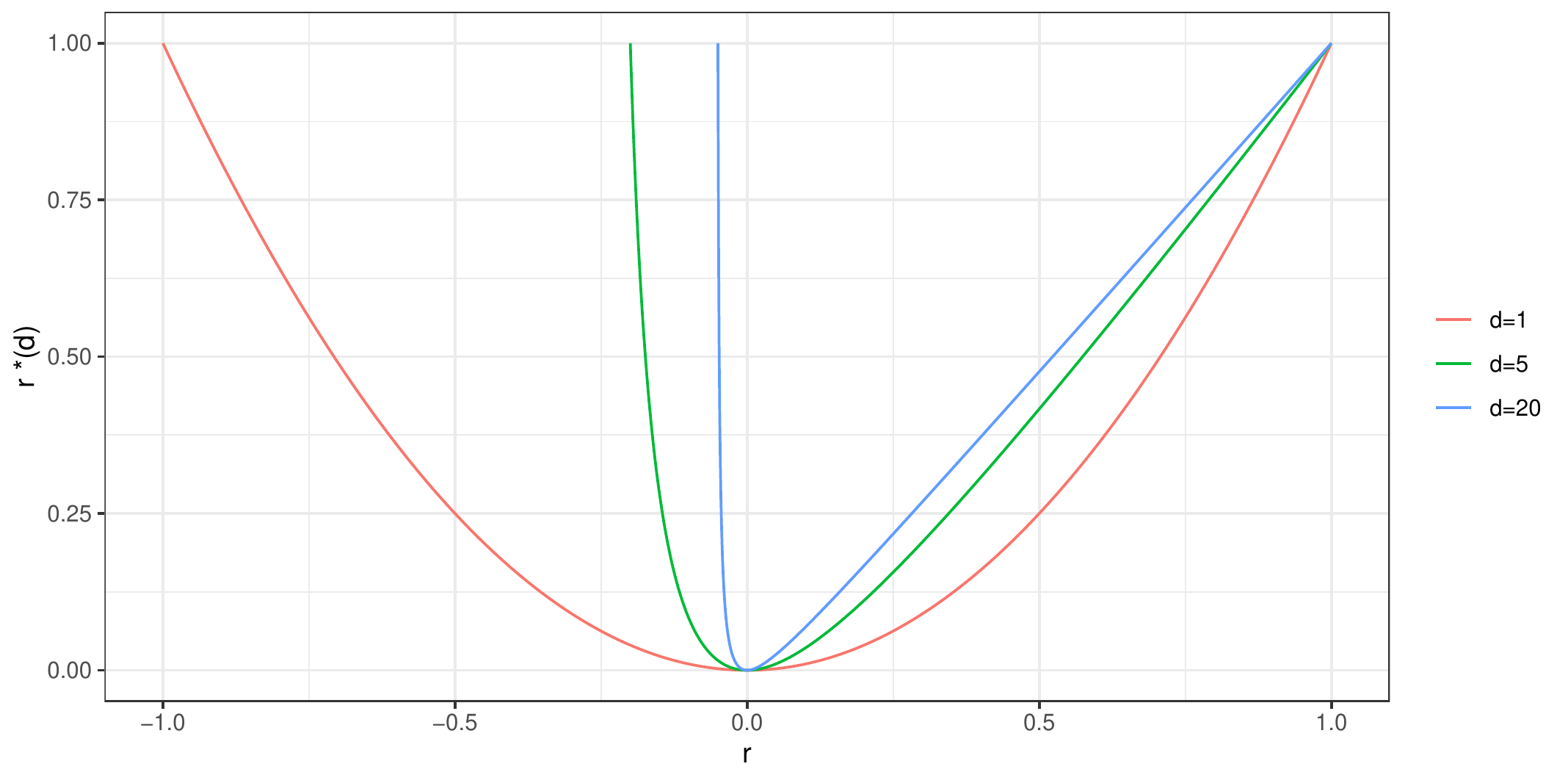}
		\caption{Plot of $r^\ast (d)$ for $d \in \{1,5,20\}$ as function of $r$.}
		\label{Fig.Gaussian.Cor}
\end{figure}	
\item (Marshall-Olkin copula) 
Consider the \emph{Marshall-Olkin} copula (see \cite[Section 6.4]{fdsempi2016}) $A_{\alpha,\beta} \in \mathcal{C}^{1+1}$ with parameters $\alpha,\beta \in [0,1]$ given by
$$ 
  A_{\alpha,\beta}(u,v):= \min \{u^{1-\alpha}v,uv^{1-\beta}\}, 
$$
and recall that $A_{\alpha,\beta} = \Pi$ if and only if $\min\{\alpha,\beta\}=0$.
After straightforward but tedious calculation, one then finds 
$$
  \psi(A_{\alpha,\beta})(s,t)
	= 
	\begin{cases}
	  \Pi(s,t) + \frac{\alpha^2}{1-2\alpha} \, \Pi(s,t) \left( 1 -  \Pi(s,t)^{\tfrac{\beta-2\alpha\beta}{\alpha}} \, M(s,t)^{\tfrac{2\alpha\beta-\beta}{\alpha}} \right),
		& \alpha \notin \{0,1/2,1\},
		\\
		A_{\beta,\beta} (s,t),
		& \alpha=1,
		\\
		\Pi(s,t) + \frac{\beta}{8} \, \Pi(s,t) \big( \log(M(s,t)) - \log(\Pi(s,t))\big), 
		& \alpha=1/2,
	\end{cases}
$$
implying that $\psi(A_{\alpha,\beta})$ is a Marshall-Olkin copula if $\min\{\alpha,\beta\}=0$ 
(in this case $\psi(A_{\alpha,\beta}) = \Pi = A_{\alpha,\beta}$) or $\alpha=1$ (in this case $\psi(A_{\alpha,\beta}) = A_{\beta,\beta}$).
Fig. \ref{Fig.MO} depicts a sample of size $n=10,000$ drawn from a Marshall-Olkin copula $A$ with parameters $\alpha=1$ and $\beta=0.3$ together with a sample of equal size from $\psi(A)$.
\begin{figure}[h!]
		\centering
		\includegraphics[width=0.65\textwidth]{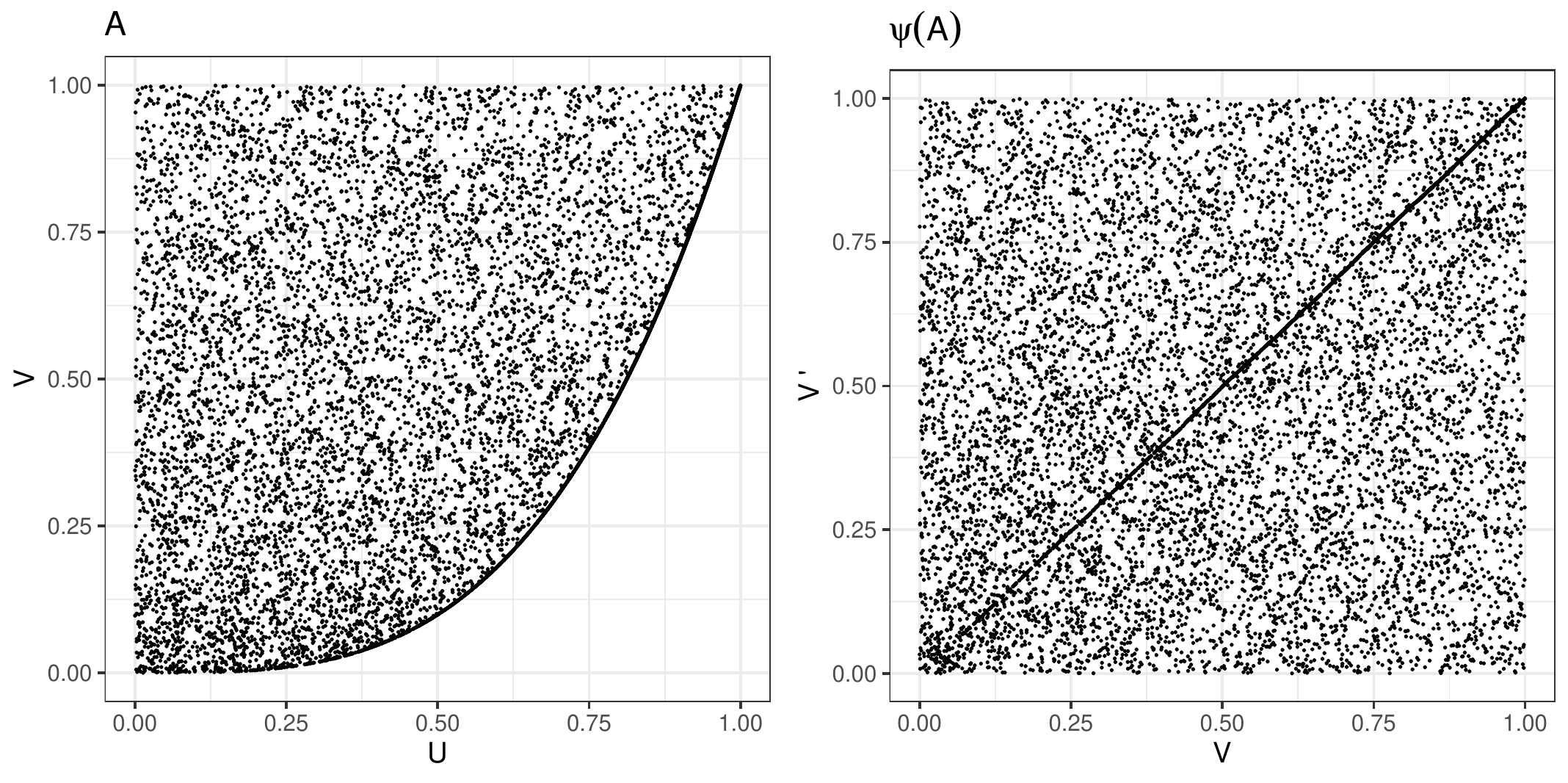}
		\caption{Sample of size $n=10,000$ drawn from a Marshall-Olkin copula $A$ with parameter $\alpha=1$ and $\beta=0.3$ (left panel), together with a sample of equal size from $\psi(A)$ (right panel).}		\label{Fig.MO}
\end{figure}
\item (Fr{\'e}chet copula) 
Consider the \emph{Fr{\'e}chet copula} (see \cite[Section 6.2]{fdsempi2016}) $A_{\alpha,\beta} \in \mathcal{C}^{1+1}$ with parameters $\alpha,\beta \in [0,1]$ such that $\alpha+\beta\leq 1$ given by
$$ 
  A_{\alpha,\beta} := \alpha \, M + (1-\alpha-\beta) \, \Pi +  \beta \, W.
$$
Then 
$ \psi(A_{\alpha,\beta})
	= A_{\alpha^\ast,\beta^\ast} $
is a Fr{\'e}chet copula with $\alpha^\ast = \alpha^2 + \beta^2$ and $\beta^\ast = 2 \alpha \beta$,
implying that the class of Fr{\'e}chet copulas is closed with respect to $\psi$.
In particular, 
$ \psi \big( A_{1/2,1/2} \big) = \psi \big( \tfrac{M+W}{2} \big) = \tfrac{M+W}{2} = A_{1/2,1/2} $.
\item (EFGM copula) 
Consider the \emph{EFGM copula} (see \cite[Section 6.3]{fdsempi2016}) $A_\alpha \in \mathcal{C}^{d+1}$ 
with parameter $\alpha \in [-1,1]$ given by 
$$  A(\uuu,v) 
	  := \Pi(\uuu,v) + \alpha \, v(1-v) \prod_{i=1}^{d} u_i(1-u_i). $$
Then 
$ \psi(A) 
    = A_{\alpha^{\ast}} $
is an EFGM copula with parameter $\alpha^\ast = \alpha^2/3^d$.
\end{enumerate}
\end{example}{}

Since $\psi(A)$ is a bivariate copula we have $\psi(A) \leq M$.
A lower bound for $\psi(A)$ can be achieved by considering the lowest possible diagonal and constructing its Bertino copula; 
we refer to \cite{frednels2003} for more information on Bertino copulas.

\begin{corollary}{} \label{phi.Diagonal}
$\psi(A)$ satisfies 
$$
  B(s,t) \leq \psi(A)(s,t) \leq M(s,t), 
	\qquad (s,t) \in \I^2,
$$
where $B$ denotes the Bertino copula with diagonal $t \to \Pi(t,t)$,
i.e., 
$$
  B (s,t) 
	= \min(s,t) - \min\{ p - p^2 \, : \, p \in [s,t]\}
	= \begin{cases}
	    \min(s,t)^2, 
		  & s+t \leq 1, 
		  \\
		  \min(s,t) - \max(s,t) + \max(s,t)^2,
		  & s+t > 1.
	  \end{cases}
$$
In particular, $\Pi(t,t) \leq \psi(A)(t,t)$ for all $t \in \I$.
\end{corollary}{}

\noindent
Notice that, although the Bertino copula $B$ with diagonal $t \to \Pi(t,t)$ serves as a lower bound of $\psi(\mathcal{C}^{d+1})$, 
it is not a member of this class as shown in Subsection \ref{Subsect.Indiff.}.

Example \ref{phi.CounterPi} illustrates that the inequality $\psi(A)(s,t) \geq \Pi(s,t)$ fails to hold for every $(s,t) \in \I^2$.

\begin{example}{}  \label{phi.CounterPi} 
According to Example \ref{phi.MO.Ex},
the copula $A := \frac{M+W}{2}$ satisfies $\psi(A) = A$,
hence there exists some $(s,t) = (1/4,3/4) \in \I^2$ such that 
$ \psi(A) (s,t) = 2/16 < 3/16 = \Pi (s,t) $.
\end{example}{}
\noindent
Due to Example \ref{phi.CounterPi}, 
the copula $\psi(A)$ fails to be stochastically increasing, left tail decreasing (LTD) and positive quadrant dependent (PQD), in general;
see \cite{sfx2023TP2,nelsen2006} for more information on these dependence properties.

The values of $\psi(A)$ outside the diagonal are bounded from above by the values along the diagonal; the result is immediate from H{\"o}lder's inequality:

\begin{corollary}{} \label{phi.Diagonal.2}
For every $A \in \mathcal{C}^{d+1}$, the inequality 
$ \psi(A) (s,t) \, \psi(A) (t,s)
	  =  \psi(A) (s,t)^2
	\leq \psi(A) (s,s) \, \psi(A) (t,t) $
holds for all $(s,t) \in \I^2$.
\end{corollary}{}

The copula $\psi(A)$ satisfies the \emph{information gain} inequality (see, e.g., \cite{fgwt2021}) along the diagonal,
i.e., the more conditioning variables are involved, the larger the value of $\psi(A)$ along the diagonal.

\begin{theorem}{} \label{phi.InformationGain}
Consider some $L \subseteq \{1,\dots,d\}$ with $1 \leq |L|\leq d-1$ and some $l \in \{1,\dots,d\} \backslash L$.
Then the inequality
$$
  \psi (A^{L \cup \{d+1\}}) (t,t)
	\leq \psi (A^{L \cup \{l\} \cup \{d+1\}}) (t,t)
$$
holds for all $t \in \I$.
In particular, the inequality
\begin{equation} \label{Eq.InformationGain}
  \psi(A^{1,d+1}) (t,t)
	\leq \psi(A^{1:2,d+1}) (t,t)
	\leq \dots
	\leq\psi(A^{1:l,d+1}) (t,t)
	\leq \dots
	\leq\psi(A) (t,t)
\end{equation}
holds for all $l \in \{1,\dots,d\}$ and all $t \in \I$.
\end{theorem}{}

\noindent
In general, the information gain inequality does not hold outside the diagonal.

\begin{example}{}
Consider the copula $A: \I^3 \to \I$ given by 
$ A(u_1,u_2,v)
	:= u_1 \, \tfrac{M+W}{2} (u_2,v) $.
Then $A^{13} = \Pi$, hence $\psi(A^{13}) = \Pi$, and, according to Example \ref{phi.MO.Ex},
$$  \psi(A)(s,t)
	= \int_{\I^2} K_{(M+W)/2} (u_2, [0,s]) \, K_{(M+W)/2} (u_2, [0,t]) \; \mathrm{d} \lambda^2 (\uuu)
	= \psi\big(\tfrac{M+W}{2}\big)(s,t)
	= \tfrac{M+W}{2} (s,t),
	\qquad (s,t) \in \I^2.
	$$
In accordance with Theorem \ref{phi.InformationGain} and Corollary \ref{phi.Diagonal}, 
we obtain
$ \psi(A^{13}) (t,t) \leq \psi(A)(t,t) $
for all $t \in \I$, however, by Example \ref{phi.CounterPi} there exists some $(s,t) = (1/4,3/4) \in \I^2$ such that $ \psi(A) (s,t) = 2/16 < 3/16 = \psi(A^{13}) (s,t) $.
\end{example}{}

Looking at \eqref{Eq.InformationGain} from the perspective of a possible dimension reduction, 
the question arises under which conditions on $A$ the information gain inequality (after a certain step) becomes an equality,
i.e., no information is added by considering additional explanatory variables.
The next result solves this question; Theorem \ref{phi.Ind.AssII} is immediate from Lemma \ref{Lemma.CI}.

\begin{theorem}{} \label{phi.Ind.AssII}
Consider a $(d+1)$-dimensional random vector $(\XXX,Y)$ with continuous marginals and connecting copula $A$.
Then, for $k \in \{1,\dots,d-1\}$, the following statements are equivalent:
\begin{enumerate}
\item[\rm (i)]
$Y$ and the random vector $(X_{k+1}, \dots, X_d)$ are conditionally independent given $(X_1, \dots, X_k)$;
\item[\rm (ii)]
$\psi(A)(t,t) = \psi(A^{1:k,d+1})(t,t)$ for all $t \in \I$.
\end{enumerate}
In addition, (a) even implies $\psi(A) = \psi(A^{1:k,d+1})$.
\end{theorem}{}

\begin{remark}{} \leavevmode 
\begin{enumerate}[1.]
\item
In view of the hierarchical feature selection performed in Subsection \ref{RDE.Subsect.FS},
Theorem \ref{phi.Ind.AssII} is of particular interest, 
as it allows to derive additional information about the underlying conditional (in-)dependence structure of the involved random variables 
when no improvement in predictability is seen after a certain step.
\item
Notice that, if $Y$ and the random vector $(X_{2}, \dots, X_d)$ are conditionally independent given $X_1$, the information gain inequality \eqref{Eq.InformationGain} becomes an equality.
This assumption is even weaker than the well-known conditional independence assumption (see, e.g., \cite{bayramoglu2014,fgwt2021}),
requiring that the random variables $X_2, \dots, X_d, Y$ are conditionally independent given $X_1$.
\end{enumerate}
\end{remark}{}

The next result shows that the copula $\psi(A)$ is invariant under a variety of measurable and bijective transformations of the first $d$ coordinates of $A$; 
this includes permutations and reflections of copulas; 
we refer to \cite{sfx2014gamma2} for more background on permutations and reflections of copulas:

\begin{corollary}{} \label{phi.Trans}
For $A \in \mathcal{C}^{d+1}$, 
consider the identity map ${\rm id}: \I \to \I$ and
some measurable bijective transformation $\zeta: \I^d \to \I^d$ such that the push-forward $(\mu_{A^{1:d}})^\zeta$ is again a copula measure.
Denoting by $A_{(\zeta,{\rm id})}$ the copula of $(\mu_{A})^{(\zeta,{\rm id})}$, 
we have $\psi(A_{(\zeta,{\rm id})}) = \psi(A)$.
\end{corollary}{}

\begin{remark}{} \label{phi.Invarianz}
In terms of a random vector $(\XXX,Y)$ with continuous marginals and connecting copula $A$, Corollary \ref{phi.Trans} implies that the transformation $\psi$ is invariant 
\begin{enumerate}[1.]
\item 
with respect to permutations of $\XXX$, and
\item
with respect to coordinatewise continuous and strictly increasing (or decreasing) transformations of $\XXX$.
\end{enumerate}
Note that $\psi$ is also invariant with respect to the linkage transformation of $(\XXX,Y)$ (considered, e.g., in \cite{lss1996, fgwt2021}) which allows to transform the random vector $(\XXX,Y)$ to a random vector $(\UUU,V)$ with uniform univariate marginals such that $U_1, \dots, U_d$ are independent.
\end{remark}{}

\begin{remark}{} \label{Rem.density.EFGM}
If $A$ is absolutely continuous with Lebesgue density $a$, 
then $\psi(A)$ is absolutely continuous as well and satisfies
$$
  \psi(A) (s,t)
	= \int_{[0,s] \times [0,t]} \int_{\{\uuu \in \I^d : a^{1:d}(\uuu)>0\}} \frac{a(\uuu, p) \, a(\uuu, q)}{a^{1:d}(\uuu)} 
				 \; \mathrm{d} \lambda^d (\uuu) \; \mathrm{d} \lambda^2 (p,q),
	\qquad (s,t) \in \I^2.
$$
For an illustration, consider the checkerboard approximation of copula $A \in \mathcal{C}^{d+1}$: 
For $N \in \mathbb{N}$ and $({\bf i},j) \in \{1,\dots,N\}^{d} \times \{1,\dots,N\}$, 
define the (hyper-)cubes 
$ T_{j}^N 
	:= \left( \frac{j-1}{N},\frac{j}{N} \right] $, 
$ S_{{\bf i}}^N 
	:= \prod_{k=1}^{d} \left( \frac{i_k-1}{N},\frac{i_k}{N}\right] $ and
$	R_{{\bf i},j}^N 
	:= S_{{\bf i}}^N \times T_{j}^N $,
and set $\sum_{{\bf i}=1}^N := \sum_{i_1=1}^N \dots \sum_{i_d=1}^N$.
Then the checkerboard approxmation $CB_N (A)$ of $A$ with resolution $N$ is given by
$$
  CB_N (A) (\uuu,v)
	:= \int_{[0,v]} \int_{[{\bf 0},\uuu]} 
				N^{d+1} \sum_{j =1}^N \sum_{{\bf i}=1}^N \mu_A (R_{{\bf i},j}^N) 
				\, \mathds{1}_{R_{{\bf i},j}^N} (\ppp,q) 
		    \; \mathrm{d} \lambda^{d} (\ppp) \mathrm{d} \lambda(q), 
$$
and since 
$(CB_N (A))^{1:d} =  CB_N (A^{1:d})$
we obtain
\begin{eqnarray*}
  \psi(CB_N (A)) (s,t)
	& = & \int_{[0,s] \times [0,t]} \int_{\{\uuu \in \I^d : \mu_{A^{1:d}} (S_{{\bf i}}^N)>0\}} 
	      N^{d+2} \sum_{j =1}^N  \sum_{k =1}^N \sum_{{\bf i} =1}^N \; 
				\frac{\mu_A (R_{{\bf i},j}^N) \, \mu_A (R_{{\bf i},k}^N)}
				                {\mu_{A^{1:d}} (S_{{\bf i}}^N)} \; 
					 \mathds{1}_{S_{{\bf i}}^N} (\uuu) \;
					 \mathds{1}_{T_{j}^N  \times T_{k}^N} (p,q)
				 \; \mathrm{d} \lambda^d (\uuu) \; \mathrm{d} \lambda^2 (p,q)
	\\
	& = & \int_{[0,s] \times [0,t]} 
	      N^{2} \sum_{j =1}^N  \sum_{k =1}^N 
				\left( \sum_{{\bf i} =1}^N \gamma_{{\bf i},j,k} \right) \mathds{1}_{T_{j}^N  \times T_{k}^N} (p,q)
        \; \mathrm{d} \lambda^2 (p,q),
\end{eqnarray*}
where
$$
  \gamma_{{\bf i},j,k}
	:= 
	\begin{cases}
	  \frac{\mu_A (S_{{\bf i}}^N \times T_{j}^N) \, \mu_A (S_{{\bf i}}^N \times T_{k}^N)}{\mu_A (S_{{\bf i}}^N \times \I)},
		& \mu_{A^{1:d}} (S_{{\bf i}}^N) > 0,
		\\
		0,
		& \textrm{else}.
  \end{cases}
$$
Thus, $\psi$ maps every ($(d+1)$-dimensional) checkerboard copula to a bivariate checkerboard copula with the same resolution.
\end{remark}{}

%%%
In the context of Markov products, it has been recognized in \cite[Theorem 5.2.10]{fdsempi2016} that, for $d=1$, the map $A \mapsto A^\top \ast A$ fails to be continuous with respect to the topology of uniform convergence and hence uniform convergence of a sequence of copulas $(A_n)_{n \in \mathbb{N}}$ to $A$ does not automatically imply uniform convergence of the sequence $(\psi(A_n))_{n \in \mathbb{N}}$ to $\psi(A)$ (see also \cite[Theorem 6]{emura2021}).

%%%%%%%%%%%%%%%%%%%%%%%%%%%%%%%%%%%%%%%%%%%%%%%%%%%%%%%%%%%%%%%%%%%%%%%%%%%%%%%%%%%%%%%%%%%%%%%%%%%%%%%%%%%%%
%%%%%%%%%%%%%%%%%%%%%%%%%%%%%%%%%%%%%%%%%%%%%%%%%%%%%%%%%%%%%%%%%%%%%%%%%%%%%%%%%%%%%%%%%%%%%%%%%%%%%%%%%%%%%
%%%%%%%%%%%%%%%%%%%%%%%%%%%%%%%%%%%%%%%%%%%%%%%%%%%%%%%%%%%%%%%%%%%%%%%%%%%%%%%%%%%%%%%%%%%%%%%%%%%%%%%%%%%%%
\section{Applications} \label{Sec.Appl}

In the present section, three applications of the dimension reduction principle induced by $\psi$ are discussed by
applying the bivariate copula $\psi(A)$ to existing measures of bivariate dependence.
Evaluating $\psi(A)$ via
\begin{enumerate}[1.]
\item Spearman's footrule leads to (the unconditional version $T$ of) the `simple measure of conditional dependence' recently introduced by \citet{chatterjee2021} (Subsection \ref{Subsec.Chatterjee}).
\item Spearman's rho leads to a distribution-free coefficient of determination $R^2$ (Subsection \ref{Subsec.COD}).
\item Gini's gamma leads to a measure of `indifference' $Q$ being able to detect whether the dependence structures of $(\XXX,Y)$ and $(\XXX,-Y)$ coincide, a kind of `lack of association' property (Subsection \ref{Subsect.Indiff.}).
\end{enumerate}{}
For either of these resulting measures, 
we investigate their ability to identify certain facets of dependence between $Y$ and $\XXX$,
show their invariance with respect to a variety of transformations of the exogenous variables $X_1, \dots, X_d$,
and verify the information gain (in-)equality.

%%%%%%%%%%%%%%%%%%%%%%%%%%%%%%%%%%%%%%%%%%%%%%%%%%%%%%%%%%%%%%%%%%%%%%%%%%%%%%%%%%%%%%%%%%%%%%%%%%%%%%%%%%%%%
%%%%%%%%%%%%%%%%%%%%%%%%%%%%%%%%%%%%%%%%%%%%%%%%%%%%%%%%%%%%%%%%%%%%%%%%%%%%%%%%%%%%%%%%%%%%%%%%%%%%%%%%%%%%%
\subsection{A measure of predictability} \label{Subsec.Chatterjee}

Considering random vectors $(\XXX,Y)$ with continuous marginals, 
the map $T$ defined in \eqref{Chatterjee.T.Intro} possesses an alternative representation in terms of the well-known dependence measure Spearman's footrule:

\begin{theorem}{} \label{zeta.Representation}
Consider a $(d+1)$-dimensional random vector $(\XXX,Y)$ with continous marginals and connecting copula $A$.
Then $T$ fulfills
$$
  T(Y,\XXX)
	= \phi(\psi(A)),
$$
where $\phi: \mathcal{C}^2 \to \mathbb{R}$ denotes Spearman's footrule given by 
$ \phi(C)
	= 6 \, \int_{\I} C(t,t) \; \mathrm{d} \lambda(t) - 2 $.
\end{theorem}{}

\noindent 
Thus, adopting the notation used in \eqref{Phi(A).Representation}, $T$ fulfills
$$
  T(Y,\XXX)
	= 6 \, \int_{\R} P (Y \leq y, Y^{\prime} \leq y) \; \mathrm{d} P^{Y}(y) - 2,
$$
with $Y$ and $Y^\prime$ sharing the same conditional distribution and being conditionally independent given $\XXX$.

The following properties of $T$ can be derived from Theorem \ref{zeta.Representation} and 
the results in Section \ref{Sect.phi} 
(Theorem \ref{phi.InDep}, Corollary \ref{phi.Trans}, Theorem \ref{phi.InformationGain}, and Theorem \ref{phi.Ind.AssII}):

\begin{corollary}{} \label{Chatterjee.T}
Consider a $(d+1)$-dimensional random vector $(\XXX,Y)$ with continuous marginals and connecting copula $A$.
Then the following properties hold:
\begin{enumerate}[(i)]
\item
$T(Y,\XXX) = 1$ if and only if $Y$ is completely dependent on $\XXX$;

\item
$T(Y,\XXX) = 0$ if and only if $Y$ and $\XXX$ are independent;
 
\item
$T(Y,\XXX)$ is invariant with respect to permutations of $X_1, \dots, X_d$;

\item
$T(Y,\XXX)$ is invariant with respect to continuous and strictly monotone transformations of $X_1, \dots, X_d$;

\item
$T$ fulfills the information gain inequality, i.e., for all $k \in \{1,\dots,d\}$,
\begin{equation*} 
  T(Y,X_1)
	\leq T(Y,(X_1,X_2))
	\leq \dots
	\leq T(Y,(X_1,\dots, X_k))
	\leq \dots
	\leq T(Y,(X_1,\dots, X_d));
\end{equation*}

\item
For $k \in \{1,\dots,d-1\}$,
the identity 
$ T(Y,(X_1,\dots, X_d))
	= T(Y,(X_1,\dots, X_k)) $ 
holds if and only if $Y$ and the random vector $(X_{k+1}, \dots, X_d)$ are conditionally independent given $(X_1,\dots,X_k)$.
\end{enumerate}
\end{corollary}{}

\noindent
Properties (i) and (ii) from Corollary \ref{Chatterjee.T} also hold for non-continuous random vectors $(\XXX,Y)$, as shown in \cite[Theorem 2.1]{chatterjee2021}.

For the copula families discussed in Example \ref{phi.GaussianII}, $T$ can be calculated explicitly.

\begin{example}{} \label{T.Ex} \leavevmode 
\begin{enumerate}[1.]
\item (Equicorrelated Gaussian copula)  
Consider the $(d+1)$-dimensional equicorrelated Gaussian copula $A_r$ with correlation parameter $r \in (-1/d,1)$.
According to Example \ref{phi.GaussianII}, $\psi(A_r)$ is a (bivariate) Gaussian copula with correlation parameter
$ r^\ast (d) = \frac{d r^2}{1+(d-1)r} $, 
and Theorem \ref{zeta.Representation} hence yields
$$
  T(Y,\XXX)
	= \phi(\psi(A_r))
	= \phi(A_{r^\ast (d)})
	= \frac{3}{\pi} \; \arcsin\left( \frac{1+r^\ast (d)}{2} \right) - 0.5
	= \frac{3}{\pi} \; \arcsin\left( \frac{1}{2} + \frac{d r^2}{2(1+(d-1)r)} \right) - 0.5.
$$
Thus, $T(Y,\XXX)=0$ if and only if $r=0$, and $T(Y,\XXX)<1$. 
\item (Marshall-Olkin copula) 
Consider the Marshall-Olkin copula $A_{\alpha,\beta}$ with parameters $\alpha=1$ and $\beta \in [0,1]$.
According to Example \ref{phi.MO.Ex}, $\psi(A)$ satisfies
$ \psi(A_{1,\beta})
	= A_{\beta,\beta}$,
and Theorem \ref{zeta.Representation} hence yields
$$
  T(Y,\XXX)
	= \phi(\psi(A_{1,\beta}))	
	= \phi(A_{\beta,\beta})	
	%= 6 \, \int_{\I} \Pi(t,t)^{1-\beta} \, t^{\beta} \; \mathrm{d} \lambda(t) - 2
	= \frac{2 \, \beta}{3-\beta}.
$$
Thus, $T(Y,\XXX)=0$ if and only if $\beta=0$, and $T(Y,\XXX)=1$ if and only if $\beta=1$. 
\item (Fr{\'e}chet copula) 
Consider the \emph{Fr{\'e}chet copula} $A_{\alpha,\beta}$ with parameter $\alpha,\beta \in [0,1]$ such that $\alpha+\beta\leq 1$.
According to Example \ref{phi.GaussianII},  
$ \psi(A_{\alpha,\beta})
	= A_{\alpha^2 + \beta^2,2\alpha\beta} $,
and Theorem \ref{zeta.Representation} hence yields
$$
  T(Y,\XXX)
	= \phi(\psi(A_{\alpha,\beta}))	
	= \phi(A_{\alpha^2 + \beta^2,2\alpha\beta})
	%= 6 \, \int_{\I} \Pi(t,t)^{1-\beta} \, t^{\beta} \; \mathrm{d} \lambda(t) - 2
	= (\alpha+\beta)^2 - 3 \, \alpha\beta
	= (\alpha-\beta)^2 + \alpha\beta.
$$
Thus, $T(Y,\XXX)=0$ if and only if $\alpha=0=\beta$, and $T(Y,\XXX)=1$ if and only if $(\alpha,\beta) \in \{\{0,1\},\{1,0\}\}$.
\item (EFGM copula) 
Consider the EFGM copula $A_\alpha$ with parameter $\alpha \in [-1,1]$.
According to Example \ref{phi.GaussianII}, $\psi(A_\alpha)$ is of EFGM type with parameter $\alpha^2/3^d$, 
and Theorem \ref{zeta.Representation} hence yields
$$
  T(Y,\XXX)
	= \phi(\psi(A_{\alpha}))	
	= \phi(A_{\alpha^2/3^d})
	= \frac{\alpha^2}{3^{d} \cdot 5}.
$$
Thus, $T(Y,\XXX)=0$ if and only if $\alpha=0$, and $T(Y,\XXX) \leq \tfrac{1}{3^{d} \cdot 5}$.
\end{enumerate}
\end{example}{}

In Subsection \ref{RDE.Subsect.FS} we apply the measure of predictability $T$ to perform a hierarchical feature
selection for analyzing the influence of a set of thermal variables on annual precipitation in a global climate data set.

%%%%%%%%%%%%%%%%%%%%%%%%%%%%%%%%%%%%%%%%%%%%%%%%%%%%%%%%%%%%%%%%%%%%%%%%%%%%%%%%%%%%%%%%%%%%%%%%%%%%%%%%%%%%%%%%%%%%%%%%
%%%%%%%%%%%%%%%%%%%%%%%%%%%%%%%%%%%%%%%%%%%%%%%%%%%%%%%%%%%%%%%%%%%%%%%%%%%%%%%%%%%%%%%%%%%%%%%%%%%%%%%%%%%%%%%%%%%%%%%%
\subsection{A nonparametric coefficient of determination} \label{Subsec.COD}

Let $(\XXX,Y)$ be a $(d+1)$-dimensional random vector such that $Y \in L^2$. 
It is well known that the variance ${\rm var} (Y)$ of $Y$ can be decomposed via 
$$
  {\rm var} (Y)
	%= {\rm var} ( E (Y \, | \, \XXX) ) + {\rm var} ( Y - E (Y \, | \, \XXX) )
	= {\rm var} ( E (Y \, | \, \XXX) ) + E ( {\rm var} (Y \, | \, \XXX) ),
$$
where ${\rm var} (E (Y \, | \, \XXX))$ equals the part of the variance of $Y$ explained by the regression function $r(\xxx) := E (Y \, | \, \XXX = \xxx)$, and
\begin{eqnarray} \label{Eq.R2.Projection}
  R^2 
	:= R^2 (Y,\XXX)
	:= \frac{{\rm var} ( E (Y \, | \, \XXX) )}{{\rm var} (Y)}
	 %= 1 - \, \frac{{\rm var} ( Y - E (Y \, | \, \XXX) )}{{\rm var} (Y)}
	 = 1 - \, \frac{E( ( Y - E (Y \, | \, \XXX) )^2)}{{\rm var} (Y)}
\end{eqnarray}
(also known as Sobol index; see \cite{gamboa2020})
denotes the proportion of the variance that is explained by the regression function $r$.

\begin{theorem}{} \label{R2.Representation.II}
Consider a $(d+1)$-dimensional random vector $(\XXX,Y)$. 
Then
$R^2$ of $Y$ given $\XXX$ fulfills
$$ R^2(Y,\XXX)
	= \varrho_P (Y,Y^\prime), 
$$
where $\varrho_P$ denotes Pearson correlation coefficient, 
and $Y$ and $Y^\prime$ are such that they share the same conditional distribution and are conditionally independent given $\XXX$.
\end{theorem}{}

Again, assume that the random vector $(\XXX,Y)$ has continuous marginals 
$F_i$ of $X_i$, $i \in \{1,\dots,d\}$, and $G$ of $Y$ and connecting copula $A$.
Then (denoting $U_i:=F_i(X_i)$, $i \in \{1,\dots,d\}$, and $V:=G(Y)$) the distribution-free $R^2(V,\UUU)$ possesses a representation in terms of the well-known dependence measure Spearman's rho; the result is immediate from Theorem \ref{R2.Representation.II}:

\begin{corollary}{} \label{R2.Representation}
Consider a $(d+1)$-dimensional random vector $(\XXX,Y)$ with continuous marginals and connecting copula $A$.
Then
$R^2(V,\UUU)$ fulfills
$$
	R^2(V,\UUU)
	= \varrho_S (\psi(A)),
$$
where $\varrho_S: \mathcal{C}^2 \to \mathbb{R}$ denotes Spearman's rank correlation coefficient (also known as Spearman's rho) given by 
$$
  \varrho_S (C)
	 = 12 \, \int_{\I^2} C(s,t) \; \mathrm{d} \lambda^2(s,t) -3
	 = 12 \, \int_{\I^2} st \; \mathrm{d} \mu_C(s,t) -3.
$$
\end{corollary}{}

\noindent 
For the bivariate case, i.e., $d=1$, 
a copula-based representation of $R^2(V,U)$ can already be found in \citet{sungur2005}.
\citet[Theorem 1]{emura2021} have further recognized that the copula correlation $R^2(V,U)$ can be expressed as Spearman's rho of the Markov product $A^\top \ast A$, which in this case coincides with $\psi(A)$.

Again, adopting the notation used in \eqref{Phi(A).Representation}, $R^2$ fulfills
$$
  R^2(V,\UUU)
	= 12 \, \int_{\R} \int_{\R} P (Y \leq y, Y^{\prime} \leq y^{\prime}) \; \mathrm{d} P^{Y}(y) \mathrm{d} P^{Y^{\prime}}(y^{\prime}) - 3,
$$
with $Y$ and $Y^\prime$ sharing the same distribution and being conditionally independent given $\XXX$.

The following properties of $R^2(V,\UUU)$ can be derived from Corollary \ref{R2.Representation} and 
the results in Section \ref{Sect.phi} (Theorem \ref{phi.InDep}, Corollary \ref{phi.Trans}, and Theorem \ref{phi.Ind.AssII}).
The information gain inequality (i.e., reducing the number of conditioning variables reduces $R^2$) does not follow from Theorem \ref{phi.InformationGain} but from Hilbert's projection theorem.

\begin{corollary}{} \label{R2.Properties}
Consider a $(d+1)$-dimensional random vector $(\XXX,Y)$ with continous marginals and connecting copula $A$.
Then the following properties hold:
\begin{enumerate}[(i)]
\item
$R^2(V,\UUU) 
= 1$ if and only if $Y$ is completely dependent on $\XXX$;

\item
$R^2(V,\UUU) 
= 0$ if and only if $E(V|\UUU=\uuu) = \frac{1}{2}$ for $\mu_{A^{1:d}}$-almost all $\uuu \in \I^d$; 
in particular, $R^2(V,\UUU) 
= 0$ whenever $Y$ and $\XXX$ are independent;
 
\item
$R^2(V,\UUU)$ 
is invariant with respect to permutations of $X_1, \dots, X_d$;

\item
$R^2(V,\UUU)$ 
is invariant with respect to continuous and strictly monotone transformations of $X_1, \dots, X_d$;

\item
$R^2(V,\UUU)$  
fulfills the \emph{information gain inequality}, 
i.e., for all $k \in \{1,\dots,d\}$,
$$
  R^2(V,U_1)
	\leq R^2(V,(U_1,U_2))
	\leq \dots 
	\leq R^2(V,(U_1,\dots,U_k))
	\leq \dots
	\leq R^2(V,(U_1,\dots,U_d));
$$
additionally,
if, for some $k \in \{1,\dots,d-1\}$, $Y$ and the random vector $(X_{k+1}, \dots, X_d)$ are conditionally independent given $(X_1,\dots,X_k)$,
then
$R^2(V,(U_1,\dots, U_d))
  = R^2(V,(U_1,\dots, U_k))$. 
\end{enumerate}
\end{corollary}{}

\noindent
For $d=1$, properties (i) and (ii) in Corollary \ref{R2.Properties} are given in \cite[Proposition 2]{emura2021}.

Notice that $R^2(Y,\XXX)$ and its distribution-free version $R^2(V,\UUU)$ may differ:

\begin{example}{} \label{Ex.Tent}Consider the case $d=1$, the uniformly distributed random variable $V \sim {\rm U}[0,1]$ and define
$$
  U := 2V \, \mathds{1}_{[0,0.5]}(V) + (2-2V) \mathds{1}_{(0.5,1]}(V) \sim {\rm U}[0,1].
$$
Then $E(V|U=u) = \frac{1}{2}$ for $\lambda$-almost all $u \in \I$ and by Corollary \ref{R2.Properties} we obtain $R^2(V,U) = 0$. 
Now, consider the random vector $(X,Y) := (U,V^2)$.
Since the map $v \mapsto v^2$ is continuous and strictly increasing on $\I$, the random vectors $(X,Y)$ and $(U,V)$ share the same copula and hence $(F(X),G(Y)) \sim (U,V)$.
In addition, 
$$
  R^2 (Y,X)
	= \frac{{\rm var} ( E (Y \, | \, X ))}{{\rm var} (Y)}
  = \frac{{\rm var} ( E (V^2 \, | \, U ))}{{\rm var} (V^2)}
	= \frac{45}{4} \, \int_{\I} \big( E (V^2 \, | \, U = u ) - E (V^2) \big)^2 \; \mathrm{d} \lambda(u)
	%= \frac{45}{4} \, \frac{1}{180}
	= \frac{1}{16},
$$
and hence 
$ R^2(V,U) = 0 < \frac{1}{16} = R^2 (Y,X) $.
\end{example}{}

For the copula families discussed in Example \ref{phi.GaussianII}, $R^2$ can be calculated explicitly.

\begin{example}{} \label{R2.Ex} \leavevmode 
\begin{enumerate}[1.]
\item (Equicorrelated Gaussian copula)  
Consider the $(d+1)$-dimensional equicorrelated Gaussian copula $A$ with correlation paramenter $r \in (-1/d,1)$.
According to Example \ref{phi.GaussianII}, $\psi(A_r)$ is a (bivariate) Gaussian copula with correlation paramenter
$ r^\ast (d) = \frac{d r^2}{1+(d-1)r} $, and Corollary \ref{R2.Representation} hence yields
$$
  R^2(V,\UUU)
	= \rho_S(\psi(A_r))
	= \rho_S(A_{r^\ast (d)})
	= \frac{6}{\pi} \; \arcsin\left( \frac{r^\ast (d)}{2} \right) 
	= \frac{6}{\pi} \; \arcsin\left( \frac{d r^2}{2(1+(d-1)r)} \right).
$$
Thus, $R^2(V,\UUU)=0$ if and only if $r=0$, and $R^2(V,\UUU)<1$. %if and only if $r^\ast (d)=1$.
\item (Marshall-Olkin copula) 
Consider the Marshall-Olkin copula $A_{\alpha,\beta}$ with parameters $\alpha=1$ and $\beta \in [0,1]$.
According to Example \ref{phi.MO.Ex}, $\psi(A)$ satisfies
$ \psi(A_{1,\beta})
	= A_{\beta,\beta}$, 
and Corollary \ref{R2.Representation} hence yields
$$
  R^2(V,\UUU)
	= \rho_S(\psi(A_{1,\beta}))	
	= \rho_S(A_{\beta,\beta})	
	= \frac{3 \, \beta}{4-\beta}.
$$
Thus, $R^2(V,\UUU)=0$ if and only if $\beta=0$, and $R^2(V,\UUU)=1$ if and only if $\beta=1$. 
\item (Fr{\'e}chet copula) 
Consider the \emph{Fr{\'e}chet copula} $A_{\alpha,\beta}$ with parameters $\alpha,\beta \in [0,1]$ such that $\alpha+\beta\leq 1$.
According to Example \ref{phi.MO.Ex},  
$ \psi(A_{\alpha,\beta})
	= A_{\alpha^2 + \beta^2,2\alpha\beta} $, 
and Corollary \ref{R2.Representation} hence yields
$$
  R^2(V,\UUU)
	= \rho_S(\psi(A_{\alpha,\beta}))	
	= \rho_S(A_{\alpha^2 + \beta^2,2\alpha\beta})
	= (\alpha - \beta)^2.  
$$
Thus, $R^2(V,\UUU)=0$ if and only if $\alpha=\beta$, and $R^2(V,\UUU)=1$ if and only if $(\alpha,\beta) \in \{\{0,1\},\{1,0\}\}$.
\item (EFGM copula) 
Consider the EFGM copula $A_\alpha$ with parameter $\alpha \in [-1,1]$.
According to Example \ref{phi.MO.Ex}, $\psi(A_\alpha)$ is of EFGM type with parameter $\alpha^2/3^d$,
and Corollary \ref{R2.Representation} hence yields
$$
  R^2(V,\UUU)
	= \rho_S(\psi(A_{\alpha}))	
	= \rho_S(A_{\alpha^2/3^d})
	= \frac{\alpha^2}{3^{d+1}}.
$$
Thus, $R^2(V,\UUU)=0$ if and only if $\alpha=0$, and $R^2(V,\UUU) \leq \tfrac{1}{3^{d+1}}$.
\end{enumerate}
\end{example}{}

The distribution-free coefficient of determination provides a benchmark for the proportion of variance that can be explained by a copula-based (regression) model. 
In Subsection \ref{RDE.Subsect.EV} we illustrate how this information can be used to judge the appropriateness of a selected copula family in the semiparametric copula-based regression model introduced by \citet{noh2013}.

%%%%%%%%%%%%%%%%%%%%%%%%%%%%%%%%%%%%%%%%%%%%%%%%%%%%%%%%%%%%%%%%%%%%%%%%%%%%%%%%%%%%%%%%%%%%%%%%%%%%%%%%%%%%%%%%%%%%%%%%
%%%%%%%%%%%%%%%%%%%%%%%%%%%%%%%%%%%%%%%%%%%%%%%%%%%%%%%%%%%%%%%%%%%%%%%%%%%%%%%%%%%%%%%%%%%%%%%%%%%%%%%%%%%%%%%%%%%%%%%%
\subsection{A measure of `indifference' or `reflection invariance'} \label{Subsect.Indiff.}

We now evaluate $\psi(A)$ via Gini's gamma which leads to a measure of `indifference' being able to detect whether the dependence structures of $(\XXX,Y)$ and $(\XXX,-Y)$ coincide (we call this property reflection invariance with respect to variable $Y$),
a kind of `lack of association' property that is fulfilled in case $\XXX$ and $Y$ are independent:
The term `lack of association' (indifference) has been introduced in an even stricter sense  by \citet{ccr1996}; more precisely, for bivariate distribution functions for which the dependence structures of $(X,Y)$, $(X,-Y)$ and even $(-X,Y)$ coincide.
Notice that the copulas of $(\XXX,Y)$ and $(\XXX,-Y)$ coincide if and only if 
$ A(\uuu,v) = A(\uuu,1) - A(\uuu,1-v) $
for all $(\uuu,v) \in \I^{d+1}$ (see, e.g., \cite{sfx2019refl}), which is equivalent to $K_A(\uuu,[0,v]) = 1- K_A(\uuu,[0,1-v])$ for $\mu_{A^{1:d}}$-almost all $(\uuu,v) \in \I^{d+1}$.

We define the map $Q$ by letting
\begin{equation} \label{Q.Gamma}
  Q(Y,\XXX)
	:= \gamma(\psi(A)),
\end{equation}
where $\gamma: \mathcal{C}^2 \to \mathbb{R}$ denotes Gini's gamma given by 
$ \gamma(C)
	= 4 \, \int_{\I} C(t,t) + C(t,1-t) \; \mathrm{d} \lambda(t) - 2 $.
The following properties of $Q(Y,\XXX)$ can be derived from \eqref{Q.Gamma} and 
the results in Section \ref{Sect.phi} (Theorem \ref{phi.InDep}, Corollary \ref{phi.Trans}).

\begin{theorem}{} \label{Q.Properties}
Consider a $(d+1)$-dimensional random vector $(\XXX,Y)$ with continous marginals and connecting copula $A$.
Then 
$$
  Q(Y,\XXX)
	= 2 \, \int_{\I} \int_{\I^d} \big( K_A(\uuu,[0,t]) - (1 - K_A(\uuu,[0,1-t])) \big)^2 \; \mathrm{d} \mu_{A^{1:d}} (\uuu) \, \mathrm{d} \lambda(t)
	\geq 0,
$$
and the following properties hold:
\begin{enumerate}[(i)]
\item
$Q(Y,\XXX) 
= 1$ if and only if $Y$ is completely dependent on $\XXX$;
\item
$Q(Y,\XXX)
= 0$ if and only if $(\XXX,Y)$ and $(\XXX,-Y)$ have the same copula;
in particular, $Q(Y,\XXX) = 0$ if $Y$ and $\XXX$ are independent;
\item
$Q(Y,\XXX)$ 
is invariant with respect to permutations of $X_1, \dots, X_d$;
\item
$Q(Y,\XXX)$ 
is invariant with respect to continuous and strictly monotone transformations of $X_1, \dots, X_d$.
\end{enumerate}
\end{theorem}{}

For the copula families discussed in Example \ref{phi.GaussianII}, $Q$ can be calculated explicitly.

\begin{example}{} \label{Q.Ex} \leavevmode 
\begin{enumerate}[1.]
\item (Equicorrelated Gaussian copula)  
Consider the $(d+1)$-dimensional equicorrelated Gaussian copula $A_r$ with correlation paramenter $r \in (-1/d,1)$.
According to Example \ref{phi.GaussianII}, $\psi(A_r)$ is a (bivariate) Gaussian copula with correlation paramenter
$ r^\ast (d) = \frac{d r^2}{1+(d-1)r} $, and Theorem \ref{Q.Properties} hence yields
$$
  Q(Y,\XXX)
	= \gamma(\psi(A_r))
	= \gamma(A_{r^\ast (d)})
	= \frac{2}{\pi} \; \arcsin\left( \frac{1+r^\ast (d)}{2} \right) - \frac{2}{\pi} \; \arcsin\left( \frac{1-r^\ast (d)}{2} \right).
$$
Thus, $Q(Y,\XXX)=0$ if and only if $r=0$, and $Q(Y,\XXX)<1$. 
\item (Marshall-Olkin copula) 
Consider the Marshall-Olkin copula $A_{\alpha,\beta}$ with parameters $\alpha=1$ and $\beta \in [0,1]$.
According to Example \ref{phi.MO.Ex}, $\psi(A)$ satisfies
$ \psi(A_{1,\beta})
	= A_{\beta,\beta}$, 
and Theorem \ref{Q.Properties} hence yields
$$
  Q(Y,\XXX)
	= \gamma(\psi(A_{1,\beta}))	
	= \gamma(A_{\beta,\beta})	
	= \frac{(4-\beta)}{(2-\beta)(3-\beta)} \; (4 - 2^{\beta}\big) - 2.
$$
Thus, $Q(Y,\XXX)=0$ if and only if $\beta=0$, and $Q(Y,\XXX)=1$ if and only if $\beta=1$. 
\item (Fr{\'e}chet copula) 
Consider the \emph{Fr{\'e}chet copula} $A_{\alpha,\beta}$ with parameters $\alpha,\beta \in [0,1]$ such that $\alpha+\beta\leq 1$.
According to Example \ref{phi.MO.Ex},  
$ \psi(A_{\alpha,\beta})
	= A_{\alpha^2 + \beta^2,2\alpha\beta} $, 
and Theorem \ref{Q.Properties} hence yields
$$
  Q(Y,\XXX)
	= \gamma(\psi(A_{\alpha,\beta}))	
	= \gamma(A_{\alpha^2 + \beta^2,2\alpha\beta})
	= (\alpha - \beta)^2.
$$
Thus, $Q(Y,\XXX)=0$ if and only if $\alpha=\beta$, and $Q(Y,\XXX)=1$ if and only if $(\alpha,\beta) \in \{\{0,1\},\{1,0\}\}$.
\item (EFGM copula) 
Consider the EFGM copula $A_\alpha$ with parameter $\alpha \in [-1,1]$.
According to Example \ref{phi.MO.Ex}, $\psi(A_\alpha)$ is of EFGM type with parameter $\alpha^2/3^d$,
and Theorem \ref{Q.Properties} hence yields
$$
  Q(Y,\XXX)
	= \gamma(\psi(A_{\alpha}))	
	= \gamma(A_{\alpha^2/3^d})
	= \frac{4 \alpha^2}{3^{d+1}\cdot 5}.
$$
Thus, $Q(Y,\XXX)=0$ if and only if $\alpha=0$, and $Q(Y,\XXX) \leq \tfrac{4}{3^{d+1} \cdot 5}$.
\end{enumerate}
\end{example}{}

Notice that the Bertino copula $B$ with diagonal $t \mapsto \Pi(t,t)$ discussed in Corollary \ref{phi.Diagonal} fulfills
$\gamma(B) = -\tfrac{1}{3} < 0$ implying $B \notin \psi(\mathcal{C}^{d+1})$.

In Subsection \ref{RDE.Subsect.RS}, we apply the introduced measure of `indifference' to detect a reflection invariant dependence structure in the fuel spray data set discussed in \cite{coblenz2020}.

%%%%%%%%%%%%%%%%%%%%%%%%%%%%%%%%%%%%%%%%%%%%%%%%%%%%%%%%%%%%%%%%%%%%%%%%%%%%%%%%%%%%%%%%%%%%%%%%%%%%%%%%%%%%%
%%%%%%%%%%%%%%%%%%%%%%%%%%%%%%%%%%%%%%%%%%%%%%%%%%%%%%%%%%%%%%%%%%%%%%%%%%%%%%%%%%%%%%%%%%%%%%%%%%%%%%%%%%%%%
%%%%%%%%%%%%%%%%%%%%%%%%%%%%%%%%%%%%%%%%%%%%%%%%%%%%%%%%%%%%%%%%%%%%%%%%%%%%%%%%%%%%%%%%%%%%%%%%%%%%%%%%%%%%%%
\section{Estimation} \label{Sect.Estimation}

We propose an estimator for $\psi(A)$ whose form is reminiscent of the empirical copula, 
but which is actually based on the graph-based estimation procedure developed by \citet{chatterjee2021}.
We show that this copula estimator is strongly consistent from which strong consistency of the plug-in estimators of Spearman's footrule, Spearman's rho and Gini's gamma can be derived.

To this end, we consider a $(d+1)$-dimensional random vector $(\XXX,Y)$
with continuous univariate marginal distribution functions $F_i$ of $X_i$, $i \in \{1,\dots,d\}$, and $G$ of $Y$ and connecting copula $A$.
Further, let $(\XXX_1,Y_1), \dots, (\XXX_n,Y_n)$ be i.i.d. copies of $(\XXX,Y)$.
Since the univariate marginals are continuous ties only occur with probability $0$.  %w.l.o.g. we can assume that there are no ties.
For each $i$, we denote by $N(i)$ the index $j$ such that $\XXX_j$ is the nearest neighbour of $\XXX_i$ with respect to the Euclidean metric on $\R^d$. 
Since there may exist several nearest neighbours of $\XXX_i$ ties are broken at random.

For $(s,t) \in \I^2$, we define 
\begin{equation} \label{Estimation.Estimate.Alternative}
  D_n (s,t) 
	:= 
	\frac{1}{n} \, \sum_{k=1}^{n} \mathds{1}_{[0,s]} (G_n^\ast(Y_k)) \mathds{1}_{[0,t]} (G_n^\ast(Y_{N(k)})),	
\end{equation}
where $G_n^\ast$ denotes a renormalized version of the empirical distribution function of $Y_1,\dots,Y_n$, i.e.,
$ G_n^\ast(y)
	= \linebreak \frac{1}{n+1} \; \sum_{k=1}^{n} \mathds{1}_{(-\infty,y]} (Y_k) $
and note that 
$$
  \psi(A) (s,t) 
	= \int_{\I^d} K_A(\uuu,[0,s]) K_A(\uuu,[0,t]) \; \mathrm{d} \mu_{A^{1:d}} (\uuu)
	= E \big( P \big( G(Y) \leq s | \XXX \big) \, P \big( G(Y) \leq t | \XXX\big) \big).
$$
By adapting the ideas developed in \cite{chatterjee2021},
we prove consistency of our estimator \eqref{Estimation.Estimate.Alternative}; 
the proof of the main Theorem \ref{Consistency.Thm2} can be found in the Appendix.

\begin{theorem}{} \label{Consistency.Thm2}
Let $(\XXX_1,Y_1), (\XXX_2,Y_2),\ldots$ be a random sample from $(\XXX,Y)$ with continuous marginals and connecting copula $A \in \mathcal{C}^{d+1}$.
Then, for all $(s,t) \in \I^2$, we have
$ \lim_{n \to \infty} D_n (s,t) = \psi(A) (s,t) $
almost surely.
\end{theorem}{}
\noindent
In particular, in the bivariate case, i.e., when $d=1$, the estimator $D_n$ is a strongly consistent estimator for the Markov product $A^\top \ast A$.

%%%%%%%%%%%%%%%%%%%%%%%%%%%%%%%%%%%%%%%%%%%%%%%%%%%%%%%%%%%%%%%%%%%%%%%%%%%%%%%%%%%%%%%%%%%%%%%%%%%%%%%%%%%%%
%%%%%%%%%%%%%%%%%%%%%%%%%%%%%%%%%%%%%%%%%%%%%%%%%%%%%%%%%%%%%%%%%%%%%%%%%%%%%%%%%%%%%%%%%%%%%%%%%%%%%%%%%%%%%
%%%%%%%%%%%%%%%%%%%%%%%%%%%%%%%%%%%%%%%%%%%%%%%%%%%%%%%%%%%%%%%%%%%%%%%%%%%%%%%%%%%%%%%%%%%%%%%%%%%%%%%%%%%%%
\subsection{Simulation study}

We illustrate the small and moderate sample performance of our estimator $D_n$ 
for two dependence structures discussed in Example \ref{phi.GaussianII}: 
(i) the equicorrelated Gaussian copula for $d=1$ with correlation parameter $\rho = 0.6$, and 
(ii) the Marshall Olkin copula with parameter vector $(\alpha,\beta) = (1,0.4)$.
We mainly restrict ourselves to the case $d=1$ in order to be able to interpret the results obtained also as estimators for the Markov product $A^\top \ast A$.

\begin{figure}[h!]
		\centering
		\includegraphics[width=1.0\textwidth]{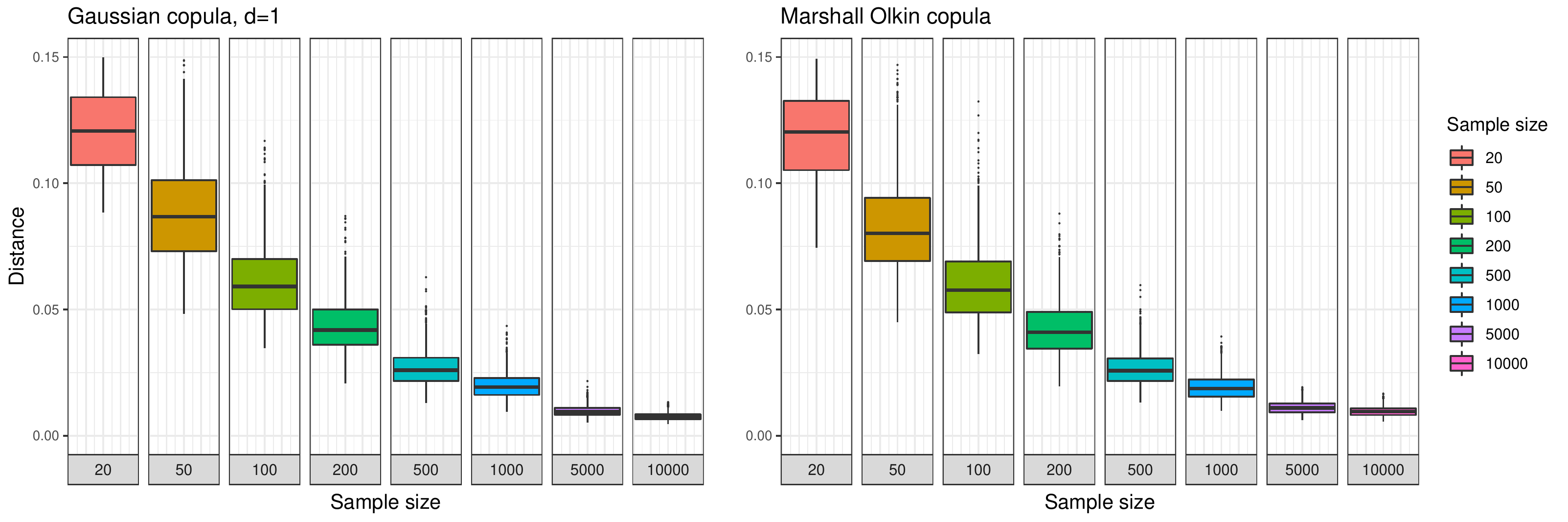}
		\caption{Boxplots summarizing the $1,000$ obtained $d_\infty$-distances between our estimate and the true copula.
		Samples of size $n$ are drawn from a Gaussian copula $A$ with parameter $\rho=0.6$ (left panel), 
		and from a Marshall-Olkin copula $A$ with parameter vector $(\alpha,\beta) = (1,0.4)$ (right panel).
		}
		\label{Fig.G1.Sim}\label{Fig.MO1.Sim}
\end{figure}

If $X$ and $Y$ have Gaussian copula $A$ with parameter $\rho = 0.6$ as connecting copula,
by Example \ref{phi.GaussianII}, the copula $\psi(A)$ is Gaussian as well with parameter $0.36$.
To test the performance of our estimator $D_n$ in this setting, 
we generated samples of size $n \in \{20; 50; 100; 200; 500; 1,000; 5,000; 10,000\}$ and calculated $D_n$.
These steps were repeated $R=1,000$ times.
Fig. \ref{Fig.G1.Sim} (left panel) depicts the $d_\infty$-distance between our estimate and the true copula evaluated on a grid of size $50$.

If $X$ and $Y$ have Marshall-Olkin copula $A$ with parameter vector $(\alpha,\beta) = (1,0.4)$ as connecting copula,
by Example \ref{phi.GaussianII}, the copula $\psi(A)$ is Marshall-Olkin with parameter vector $(0.4,0.4)$.
To test the performance of our estimator $D_n$ in this setting, 
we generated samples of size $n \in \{20; 50; 100; 200; 500; 1,000; 5,000; 10,000\}$ and calculated $D_n$.
These steps were repeated $R=1,000$ times.
Fig. \ref{Fig.MO1.Sim} (right panel) depicts the $d_\infty$-distance between our estimate and the true copula evaluated on a grid of size $50$.

As can be seen from Fig. \ref{Fig.G1.Sim}, 
the copula estimate converges rather fast to the true copula.

%%%%%%%%%%%%%%%%%%%%%%%%%%%%%%%%%%%%%%%%%%%%%%%%%%%%%%%%%%%%%%%%%%%%%%%%%%%%%%%%%%%%%%%%%%%%%%%%%%%%%%%%%%%%%
%%%%%%%%%%%%%%%%%%%%%%%%%%%%%%%%%%%%%%%%%%%%%%%%%%%%%%%%%%%%%%%%%%%%%%%%%%%%%%%%%%%%%%%%%%%%%%%%%%%%%%%%%%%%%
%%%%%%%%%%%%%%%%%%%%%%%%%%%%%%%%%%%%%%%%%%%%%%%%%%%%%%%%%%%%%%%%%%%%%%%%%%%%%%%%%%%%%%%%%%%%%%%%%%%%%%%%%%%%%
\subsection{Applications}

Finally, we plug our consistent estimator \eqref{Estimation.Estimate.Alternative} into the functionals Spearman's footrule, Spearman's rho and Gini's gamma (discussed in Section \ref{Sec.Appl}), 
which then leads to consistent estimators for the maps $T$, $R^2$ and $Q$.

%%%%%%%%%%%%%%%%%%%%%%%%%%%%%%%%%%%%%%%%%%%%%%%%%%%%%%%%%%%%%%%%%%%%%%%%%%%%%%%%%%%%%%%%%%%%%%%%%%%%%%%%%%%%%
%%%%%%%%%%%%%%%%%%%%%%%%%%%%%%%%%%%%%%%%%%%%%%%%%%%%%%%%%%%%%%%%%%%%%%%%%%%%%%%%%%%%%%%%%%%%%%%%%%%%%%%%%%%%%
\subsubsection*{A measure of predictability}

\noindent
As estimator for $T$ Akzadia \& Chatterjee \cite{chatterjee2021} propose to use the statistic
\begin{eqnarray*}
  T_n(Y,\XXX)
	& = & \frac{\sum_{i=1}^{n} (n \, \min\{R_i,R_{N(i)}\}-L_i^2)}{\sum_{i=1}^{n} L_i \,(n - L_i)},
\end{eqnarray*}
where $R_i$ denotes the rank of $Y_i$ among $Y_1, \dots, Y_n$, i.e., the number of $j$ such that $Y_j \leq Y_i$, 
and $L_i$ denotes the number of $j$ such that $Y_j \geq Y_i$ (see also \cite{bernoulli2021}).
Straightforward calculation yields
\begin{eqnarray*}
  T_n(Y,\XXX)
	& = & 1 - \frac{3}{n^2-1} \sum_{i=1}^{n} |R_i - R_{N(i)}| + \frac{3}{n^2-1} \left( \sum_{i=1}^{n} R_{N(i)} + \sum_{i=1}^{n} R_{i} - n(n+1) \right).
\end{eqnarray*}
In view of Theorem \ref{zeta.Representation},
it is worth to note that $T_n$ equals
$$
  T_n(Y,\XXX)
	= \frac{n}{n-1} \; \phi(D_n) - \frac{1}{n-1},
$$
with $\phi(D_n)$ being the plug-in estimator of Spearman's footrule (see \cite{sfx2019spearman,genebg2010,nelsen2006}) and $D_n$ being the consistent estimator for $\psi(A)$ given in \eqref{Estimation.Estimate.Alternative}.
Notice that there may exist more than one index $i$ such that $\XXX_j$ is a nearest neighbour of $\XXX_i$ implying that 
$\sum_{i=1}^{n} R_{N(i)}$ may fail to equal $\frac{n(n+1)}{2}$.

It has been proven in \cite[Theorem 2.2]{chatterjee2021} that $T_n(Y,\XXX)$ is a strongly consistent estimator for $T(Y,\XXX)$.
In \cite{shi2021normal} the authors showed asymptotic normality of $\sqrt{n} T_n$ under independence and for some regularity conditions; 
for a comprehensive summary of properties for $T_n$ we refer to \citet{bernoulli2021} and the references therein.

%%%%%%%%%%%%%%%%%%%%%%%%%%%%%%%%%%%%%%%%%%%%%%%%%%%%%%%%%%%%%%%%%%%%%%%%%%%%%%%%%%%%%%%%%%%%%%%%%%%%%%%%%%%%%
%%%%%%%%%%%%%%%%%%%%%%%%%%%%%%%%%%%%%%%%%%%%%%%%%%%%%%%%%%%%%%%%%%%%%%%%%%%%%%%%%%%%%%%%%%%%%%%%%%%%%%%%%%%%%
\subsubsection*{A nonparametric coefficient of determination}

\noindent
Motivated by Corollary \ref{R2.Representation},
as estimator for $R^2(V,\UUU)$ we propose to use the statistic
$$
  R^2_n (V,\UUU)
	:= \frac{12}{n(n+1)^2} \sum_{i=1}^{n} R_i \, R_{N(i)} - 3 - \frac{12}{n(n+1)} \left( \sum_{i=1}^{n} R_{N(i)} + \sum_{i=1}^{n} R_{i} - n(n+1) \right), 
$$
which equals the plug-in estimator $\rho_S(D_n)$ of Spearman's rho (see \cite{nelsen2006}) with $D_n$ being the consistent estimator for $\psi(A)$ given in \eqref{Estimation.Estimate.Alternative}.
By Theorem \ref{Consistency.Thm2}, the estimator $R^2_n(V,\UUU)$ is a strongly consistent estimator for $R^2(V,\UUU)$.
\\
For the bivariate case, i.e., $d=1$, \citet{gamboa2020} introduced an estimator for $R^2(Y,X)$ based on Pearson's correlation coefficient using the technique developed by \citet{chatterjee2020}.

%%%%%%%%%%%%%%%%%%%%%%%%%%%%%%%%%%%%%%%%%%%%%%%%%%%%%%%%%%%%%%%%%%%%%%%%%%%%%%%%%%%%%%%%%%%%%%%%%%%%%%%%%%%%%
%%%%%%%%%%%%%%%%%%%%%%%%%%%%%%%%%%%%%%%%%%%%%%%%%%%%%%%%%%%%%%%%%%%%%%%%%%%%%%%%%%%%%%%%%%%%%%%%%%%%%%%%%%%%%
\subsubsection*{A measure of \enquote{indifference} or \enquote{reflection invariance}}

\noindent
Motivated by \eqref{Q.Gamma},
as estimator for $Q(Y,\XXX)$ we propose to use the statistic
$$
  Q_n (Y,\XXX)
	:= \frac{2}{n(n+1)} \left( \sum_{i=1}^{n} |R_i + R_{N(i)} - (n+1)| - \sum_{i=1}^{n} |R_i - R_{N(i)}| \right) + \frac{4}{n(n+1)} \left( \sum_{i=1}^{n} R_{N(i)} + \sum_{i=1}^{n} R_{i} - n(n+1) \right), 
$$
which equals the plug-in estimator $\gamma_S(D_n)$ of Gini's gamma (see \cite{nelsen2006}) with $D_n$ being the consistent estimator for $\psi(A)$ given in \eqref{Estimation.Estimate.Alternative}.
By Theorem \ref{Consistency.Thm2}, the estimator $Q_n$ is a strongly consistent estimator for $Q$.

%%%%%%%%%%%%%%%%%%%%%%%%%%%%%%%%%%%%%%%%%%%%%%%%%%%%%%%%%%%%%%%%%%%%%%%%%%%%%%%%%%%%%%%%%%%%%%%%%%%%%%%%%%%%%
%%%%%%%%%%%%%%%%%%%%%%%%%%%%%%%%%%%%%%%%%%%%%%%%%%%%%%%%%%%%%%%%%%%%%%%%%%%%%%%%%%%%%%%%%%%%%%%%%%%%%%%%%%%%%
%%%%%%%%%%%%%%%%%%%%%%%%%%%%%%%%%%%%%%%%%%%%%%%%%%%%%%%%%%%%%%%%%%%%%%%%%%%%%%%%%%%%%%%%%%%%%%%%%%%%%%%%%%%%%
\section{Real Data example} \label{Sec.DataEx}

Finally, we illustrate the potential and importance of the functionals discussed in Section \ref{Sec.Appl} by analyzing several real data examples.
In Subsection \ref{RDE.Subsect.FS}, by applying the coefficient $T$ we perform a feature selection for a data set of bioclimatic variables and at the same time determine what proportion of the variance the selected variables are capable of explaining.
In Subsection \ref{RDE.Subsect.EV}, we then point out how the distribution-free coefficient of determination can be used for judging the appropriateness of a selected copula-based (regression) model, and to which extent it can assist in choosing the best copula family out of a set of candidate families.
We conclude this section by demonstrating how the measure of `indifference' can help identifying reflection invariant dependence structures in a fuel spray data set (Subsection \ref{RDE.Subsect.RS}).

%%%%%%%%%%%%%%%%%%%%%%%%%%%%%%%%%%%%%%%%%%%%%%%%%%%%%%%%%%%%%%%%%%%%%%%%%%%%%%%%%%%%%%%%%%%%%%%%%%%%%%%%%%%%%
\subsection{Analysis of global climate data: Feature selection} \label{RDE.Subsect.FS}

We consider a data set of bioclimatic variables for $n=1862$ locations homogeneously distributed over the global landmass from CHELSEA (\cite{karger2017, karger2018}) and want to analyze the influence of a set of thermal variables on \emph{Annual Precipitation} (AP).

For this purpose, by applying the coefficient $T$ we perform a hierarchical feature selection and identify those variables that best predict AP (= variable Y).
Fig. \ref{Fig.Bio1} depicts the order of the hierarchically selected variables based on the estimated value for $T$.
There, the value in line $k$ indicates the estimated value for $T(Y,(X_1,\dots,X_k))$ where $X_1,\dots,X_k$ are the variables in lines $1$ to $k$. As an additional feature, the last column in Fig. \ref{Fig.Bio1} contains the estimated values for the distribution-free coefficient of determination $R^2$ in this model
and hence provides a benchmark for the proportion of variance that can be explained by a copula-based model.

The increase in the estimated values for $T$ and $R^2$ with increasing number of explanatory variables is in accordance with the information gain inequalities discussed in 
Corollary \ref{Chatterjee.T} and Corollary \ref{R2.Properties}.

\begin{table}[h!]
\caption{Results of the hierarchical feature selection based on the coefficient $T$ to identify those variables that best predict AP; the last column contains the estimated value for $R^2$ in this model.}
\label{Fig.Bio1}
\centering
\begin{tabular}{c|lcc}
Position 
& Variables 
& Estimate for $T$
& Estimate for $R^2$
\\ 
\hline
1 &
\emph{Min Temperature of Coldest Month} (MTCM) 
& 0.463
& 0.628
\\
2 &
\emph{Mean Diurnal Range} (MDR) 
& 0.666
& 0.821
\\
3 &
\emph{Max Temperature of Warmest Month} (MTWM) 
& 0.722
& 0.876
\\
4 &
\emph{Mean Temperature of Wettest Quarter} (MTWeQ) 
& 0.743
& 0.891
\\
5 &
\emph{Temperature Seasonality} (TS) 
& 0.754
& 0.896
\\
6 &
\emph{Mean Temperature of Driest Quarter} (MTDQ) 
& 0.762
& 0.901
\\
7 &
\emph{Mean Temperature of Warmest Quarter} (MTWaQ)) 
& 0.764
& 0.905
\\
8 &
\emph{Mean Temperature of Coldest Quarter} (MTCQ) 
& 0.763
& 0.903
\\
9 &
\emph{Isothermality} (IT) 
& 0.766
& 0.902
\\
10 &
\emph{Annual Mean Temperature} (AMT)
& 0.765
& 0.905
\\
11 &
\emph{Temperature Annual Range} (TAR) 
& 0.764
& 0.906
\\
\end{tabular}
\end{table}

\begin{figure}[h!]
		\centering
		\includegraphics[width=0.8\textwidth]{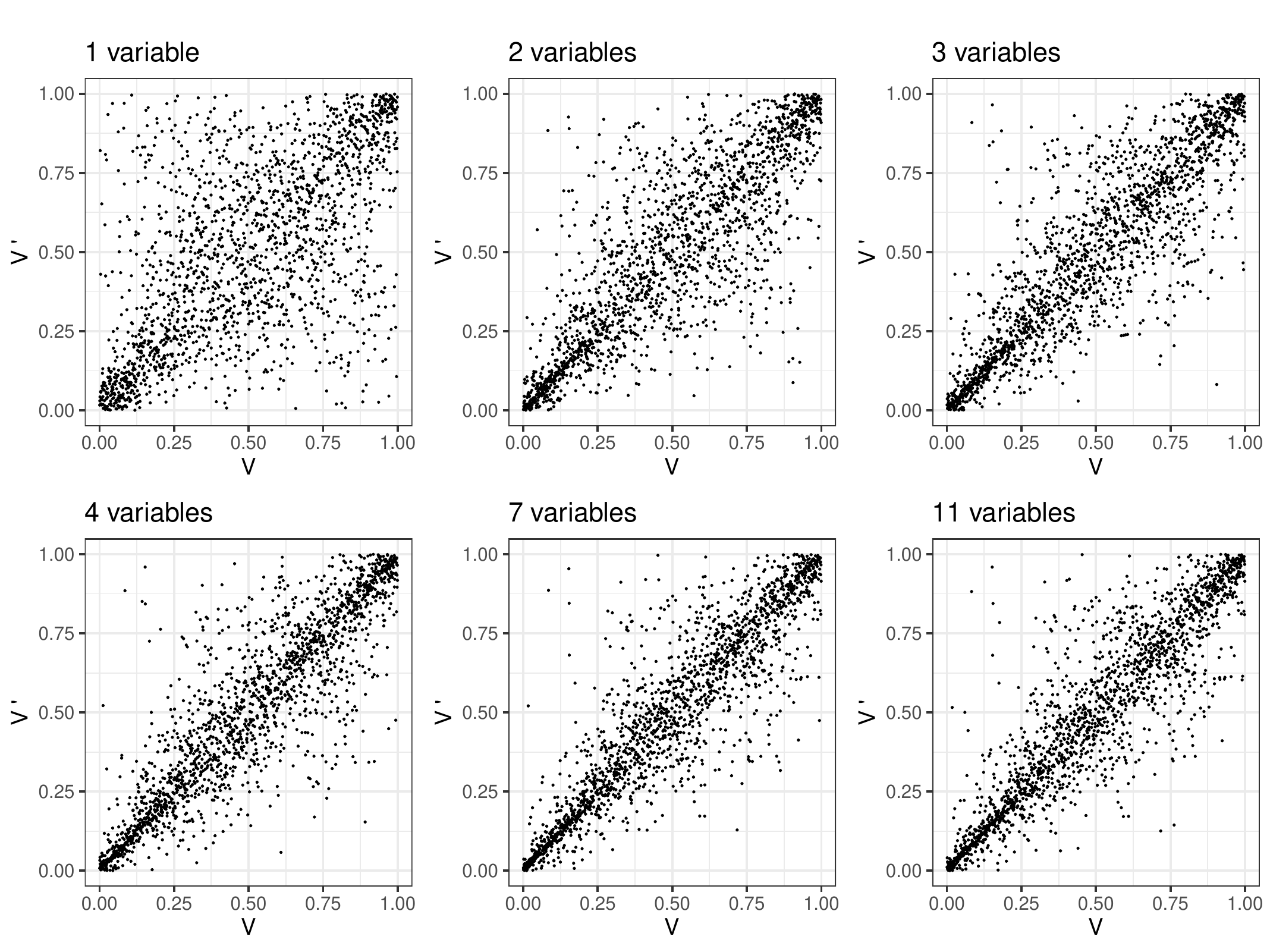}
		\caption{Sample plots of the estimate for $\psi(A)$ when the first $1$,$2$,$3$,$4$,$7$ and $11$ variables from Table \ref{Fig.Bio1} are used to predict AP.}
		\label{Fig.Bio2}
\end{figure}

Table \ref{Fig.Bio1} indicates that it is sufficient to use only a selected number of variables (here, e.g., 6) to build a model, as there is no improvement in the estimate for $T$ above this threshold. 
In view of Corollary \ref{Chatterjee.T} (vi), this also means that the variables $Y$ and $X_7, \dots, X_{11}$ can be considered as conditionally independent given the variables $X_1, \dots, X_6$.
It is further interesting to see in Fig. \ref{Fig.Bio2} that beyond a certain number of variables involved, not only do $T$ and $R^2$ remain almost constant, but also the underlying dependence structure for measuring $T$ and $R^2$ no longer changes.

%%%%%%%%%%%%%%%%%%%%%%%%%%%%%%%%%%%%%%%%%%%%%%%%%%%%%%%%%%%%%%%%%%%%%%%%%%%%%%%%%%%%%%%%%%%%%%%%%%%%%%%%%%%%%
\subsection{Copula-based regression: Explained variance} \label{RDE.Subsect.EV}

\begin{figure}[h!]
		\centering
		\includegraphics[width=0.6\textwidth]{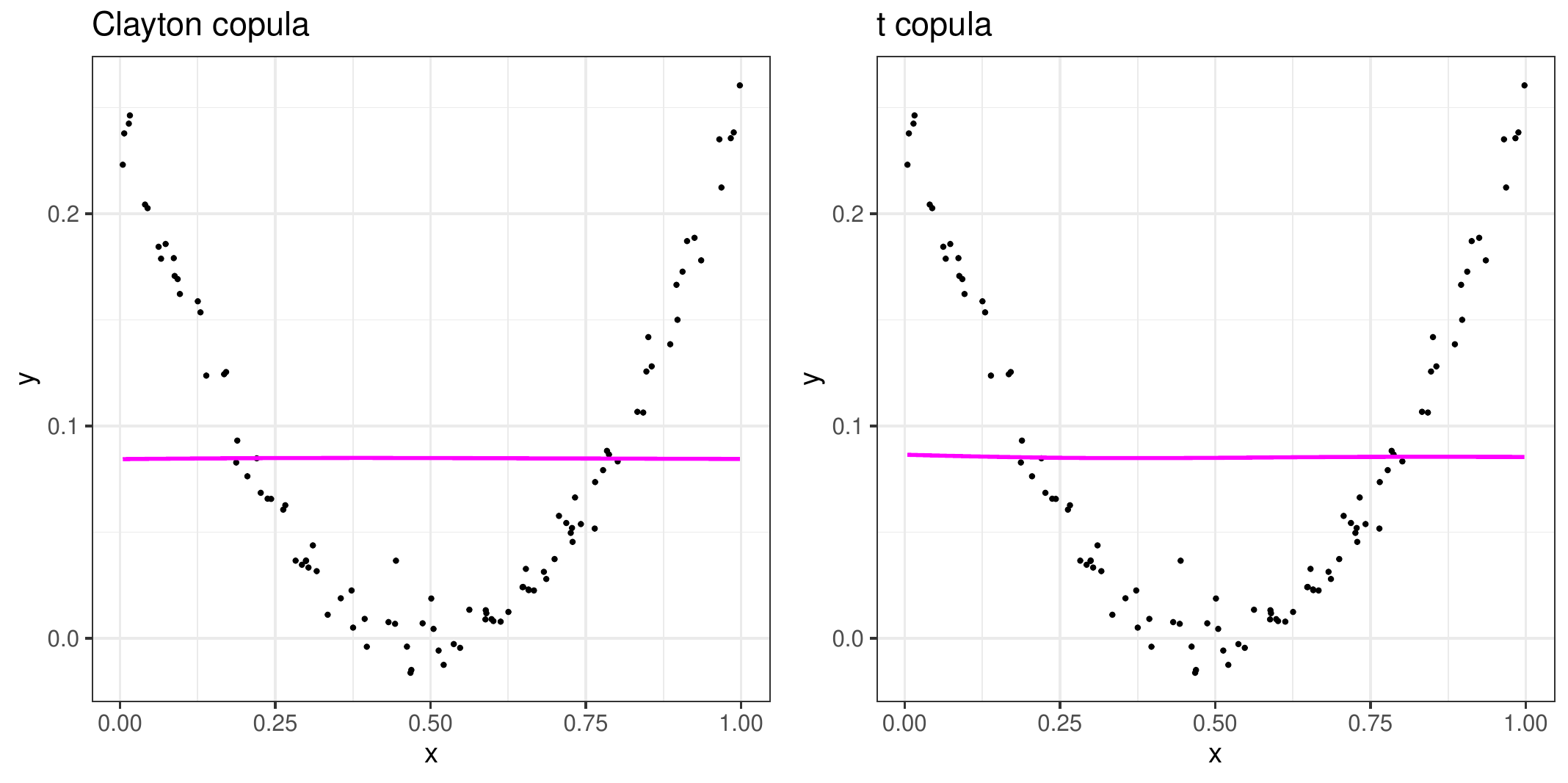}
		\caption{Semiparametric copula-based regression estimates of model \eqref{Ex.Dette} when estimating from the Clayton family of copulas (left panel) and the family of $t$ copulas (right panel).}
		\label{Fig.CopReg}
\end{figure}

We now illustrates how the distribution-free coefficient of determination discussed in Subsection \ref{Subsec.COD} can be used for judging the appropriateness of a selected copula-based regression model. 
To this end, let us consider the semiparametric copula-based regression estimator introduced by \citet{noh2013}.
For a random vector $(\XXX,Y)$ with continuous univariate marginals and connecting copula $A$, the authors showed that the mean regression function
$ r(\xxx) = E (Y \, | \, \XXX = \xxx)$ can be written as 
$$
  r(\xxx) = \int_{\R} \frac{a(F_1(x_1), \dots, F_d(x_d), G(y))}{a^{1:d}(F_1(x_1), \dots, F_d(x_d))} \, \mathrm{d} P^{Y}(y),
$$
where $a$ denotes the density of copula $A$ (also see \cite{dette2014}). 
Concerning the estimation of $r$, \citet{noh2013} suggest a semiparametric approach in which the marginal distribution functions are estimated nonparametrically and the copula is estimated parametrically from a given copula family.
The authors showed that this regression estimator is asymptotically normal if the parametric copula has been selected correctly.
However, as pointed out by \citet{dette2014} \enquote{the quality of the estimate under misspecification of the parametric copula depends heavily on the specific structure of the unknown regression function}.
\citet{dette2014} underpinned their statement by considering the simple univariate (and nonmonotone) regression model 
\begin{equation} \label{Ex.Dette}
  Y_i = (X_i - 0.5)^2 + \sigma \, \varepsilon_i, 
\end{equation}
with $X_i$ being uniformly distributed on $[0,1]$ and $\varepsilon_i$ being normally distributed with mean $0$ and variance $\sigma^2=0.01$, $i \in \{1,\dots,n\}$. 
They showed that no copula from standard copula classes \enquote{reproduces the structure of the regression function in the resulting estimate} which is due to the fact that \enquote{none of the available parametric copula models} for the vector $(X,Y)$ \enquote{yields a nonmonotone regression function}.
Indeed, although the estimated value for the distribution-free coefficient of determination $R^2$ equals $0.964$, indicating the existence of a suitable copula-based model, the explained variance when estimating from the family of $t$ copulas and the family of Clayton copulas is less than $0.01$ ($t$: 0.0046, Clayton: 0.0001) as illustrated in Fig. \ref{Fig.CopReg}.

In a second step, we illustrate the extent to which the distribution-free coefficient of determination $R^2$ can assist in choosing the best copula family out of a set of candidate families.
Therefore, let us consider the data set {\sf faithful} provided in the R package {\sf datasets}. 
The data set contains $n=272$ observations of the waiting times between eruptions (variable {\sf waiting}) and the duration of the eruption (variable {\sf eruptions}) for the Old Faithful geyser in Yellowstone National Park, Wyoming, USA.
By applying the above-mentioned semiparametric copula-based regression method, we want to estimate a regression function describing the duration of the eruptions by the waiting times between eruptions, choosing from the following three parametric copula families: `Gaussian', `Clayton' and `Joe'.
\begin{figure}[h!]
		\centering
		\includegraphics[width=0.7\textwidth]{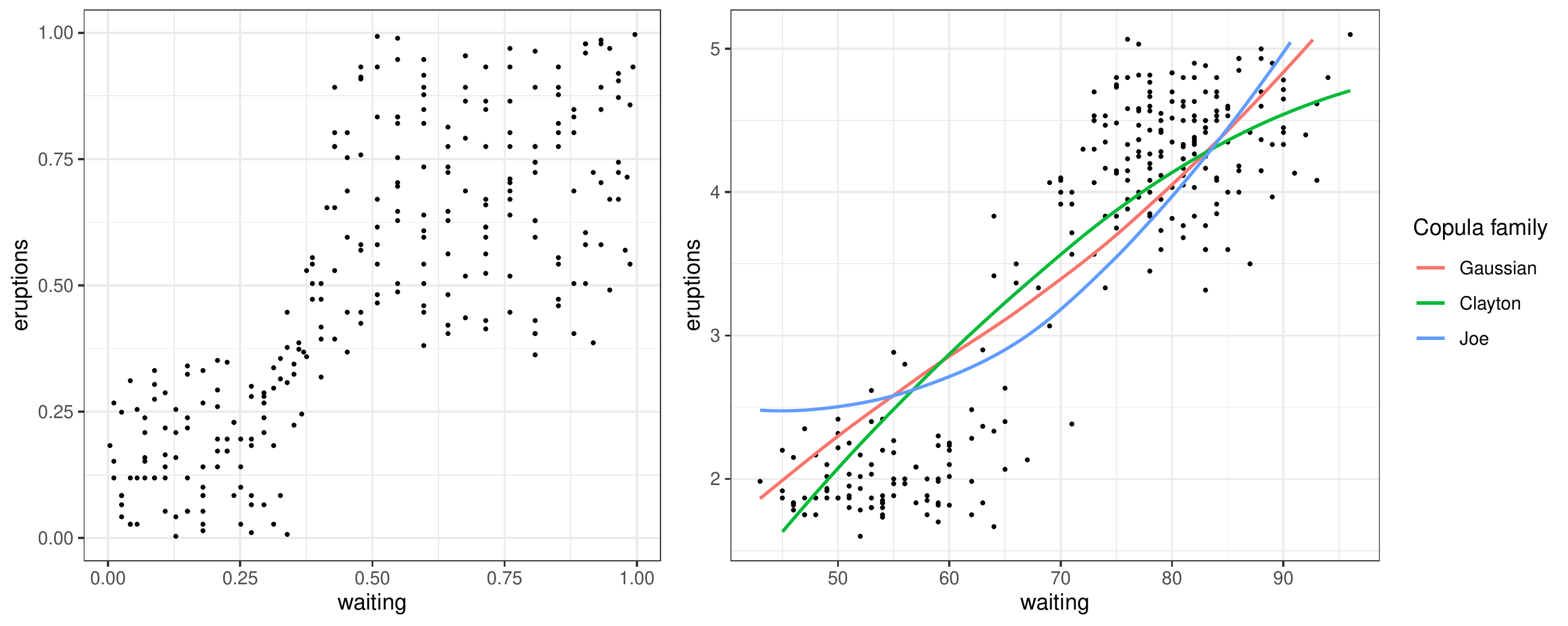}
		\caption{Pseudo-observations of data set {\sf faithful} (left panel); semiparametric copula-based regression estimates of data set {\sf faithful} when estimating from the Gaussian, Clayton, and Joe family of copulas (right panel).}
		\label{Fig.CopRegII}
\end{figure}	
Calculating the explained variance of the different copula models (Gaussian: $0.582$, Clayton: $0.633$, Joe: $0.494$) indicates that the Clayton family-based regression model performs better than its two competitors; Fig. \ref{Fig.CopRegII} depicts the regression estimates of data set {\sf faithful} when estimating from the Gaussian, Clayton, and Joe family of copulas.
However, comparing these values with the estimated value $0.715$ for the distribution-free coefficient of determination, which serves as a benchmark for the proportion of variance that can be explained by a copula-based regression model, it becomes apparent that the Clayton model may not be the optimal choice for this data set.

%%%%%%%%%%%%%%%%%%%%%%%%%%%%%%%%%%%%%%%%%%%%%%%%%%%%%%%%%%%%%%%%%%%%%%%%%%%%%%%%%%%%%%%%%%%%%%%%%%%%%%%%%%%%%
\subsection{Detecting reflection invariant dependence structures} \label{RDE.Subsect.RS}

We finally demonstrate how the measure of `indifference' introduced in Subsection \ref{Subsect.Indiff.} can be used to detect reflection invariant dependence structures.

\begin{figure}[h!]
		\centering
		\includegraphics[width=0.65\textwidth]{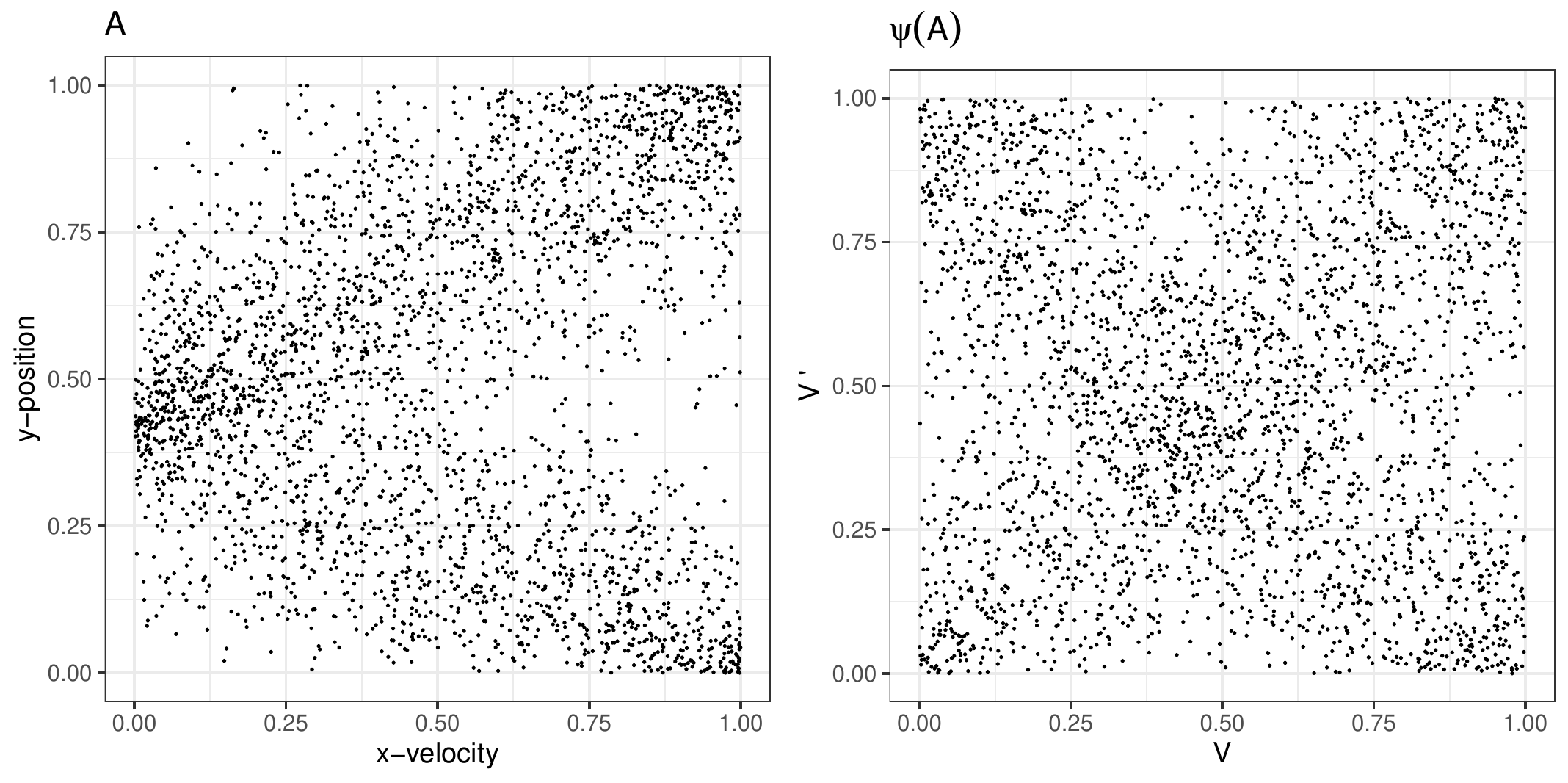}
		\caption{Pseudo-observations of the variables {\sf x-velocity} and {\sf y-position} (with connecting copula $A$) from the fuel spray data set discussed in \cite{coblenz2020} (left panel), together with a sample plot of the estimate for $\psi(A)$ (right panel).}
		\label{Fig.Gini}
\end{figure}

For this purpose, we consider the fuel spray data set discussed in \cite{coblenz2020}, which describes the behaviour of the fuel spray droplets for a specific jet engine operating condition.
According to \cite{coblenz2020}, in \enquote{jet engines the fuel is typically injected by so-called prefilming airblast atomizers} 
(we refer to \cite{coblenz2020} for more information on the physical process),
and the behaviour of the droplets is modeled using the variables {\sf drop size}, {\sf x-position}, {\sf y-position}, {\sf x-velocity}, and {\sf y-velocity}. 
Interestingly, the dependence structure of variables {\sf x-velocity} and {\sf y-position}, as depicted in Fig. \ref{Fig.Gini} (left panel), exhibits a reflection invariant dependence structure with respect to the variable {\sf y-position},
i.e., for a given velocity of the droplets along the x-coordinate, the droplets distribute symmetrically along the y-coordinate, 
a behaviour that is highly desirable for the combustion process.

The reflection invariance of the dependence structure $A$ between {\sf x-velocity} and {\sf y-position} can be easily detected from the estimated values for $T$ and $Q$. 
Since $Q_n = 0.0021$ and $T_n = 0.0589$, one may conclude that the dependence structure is reflection invariant with respect to {\sf y-position}, but the two variables are not independent; 
the latter can also be deduced from Theorem \ref{phi.InDep} and the fact that the estimate for $\psi(A)$ does not coincide with the independence copula (as seen in the right panel of Fig. \ref{Fig.Gini}).

%%%%%%%%%%%%%%%%%%%%%%%%%%%%%%%%%%%%%%%%%%%%%%%%%%%%%%%%%%%%%%%%%%%%%%%%%%%%%%%%%%%%%%%%%%%%%%%%%%%%%%%%%%%%%%%%%%
%%%%%%%%%%%%%%%%%%%%%%%%%%%%%%%%%%%%%%%%%%%%%%%%%%%%%%%%%%%%%%%%%%%%%%%%%%%%%%%%%%%%%%%%%%%%%%%%%%%%%%%%%%%%%%%%%%
%%%%%%%%%%%%%%%%%%%%%%%%%%%%%%%%%%%%%%%%%%%%%%%%%%%%%%%%%%%%%%%%%%%%%%%%%%%%%%%%%%%%%%%%%%%%%%%%%%%%%%%%%%%%%%%%%%

\section*{Acknowledgement}
The author gratefully acknowledges the support of the Austrian Science Fund (FWF) project P 36155-N 
\enquote{ReDim: Quantifying Dependence via Dimension Reduction}
and the support of the WISS 2025 project 
'IDA-lab Salzburg' (20204-WISS/225/197-2019 and 20102-F1901166-KZP).

%%%%%%%%%%%%%%%%%%%%%%%%%%%%%%%%%%%%%%%%%%%%%%%%%%%%%%%%%%%%%%%%%%%%%%%%%%%%%%%%%%%%%%%%%%%%%%%%%%%%%%%%%%%%%%%%%%
%%%%%%%%%%%%%%%%%%%%%%%%%%%%%%%%%%%%%%%%%%%%%%%%%%%%%%%%%%%%%%%%%%%%%%%%%%%%%%%%%%%%%%%%%%%%%%%%%%%%%%%%%%%%%%%%%%
%%%%%%%%%%%%%%%%%%%%%%%%%%%%%%%%%%%%%%%%%%%%%%%%%%%%%%%%%%%%%%%%%%%%%%%%%%%%%%%%%%%%%%%%%%%%%%%%%%%%%%%%%%%%%%%%%%
%\bibliography{biblio}

\begin{thebibliography}{}

\bibitem[\protect\citeauthoryear{Auddy, Deb, and Nandy}{Auddy
  et~al.}{2021}]{auddy2021}
A. Auddy, N.~Deb, S.~Nandy,
\newblock Exact detection thresholds for {C}hatterjee's correlation,
\newblock Available at \url{https://arxiv.org/abs/2104.15140v1} (2021).

\bibitem[\protect\citeauthoryear{Azadkia and Chatterjee}{Azadkia and
  Chatterjee}{2021}]{chatterjee2021}
M. Azadkia, S.~Chatterjee,
\newblock A simple measure of conditional dependence,
\newblock {Ann. Stat. 49\/} (2021) 3070--3102.

\bibitem[\protect\citeauthoryear{Bayramoglu}{Bayramoglu}{2010}]{bayramoglu2014}
I. Bayramoglu,
\newblock On conditionally independent random variables, copula and order statistics,
\newblock {Comm. Statist. Theory Meth.\/}~{43\/} (2010) 2105--2117.

\bibitem[\protect\citeauthoryear{Bickel}{Bickel}{2022}]{bickel2022}
P. Bickel,  
\newblock Measures of independence and functional dependence,
\newblock Available at \url{https://arxiv.org/abs/2206.13663v1} (2022).

\bibitem[\protect\citeauthoryear{B{\"o}ttcher, Keller-Ressel, and
  Schilling}{B{\"o}ttcher et~al.}{2019}]{bottcher2019}
B. B{\"o}ttcher, M.~Keller-Ressel, R.~Schilling,
\newblock Distance multivariance: New dependence measures for random vectors,
\newblock {Ann. Statist.\/}~{47\/} (2019) 2757--2789.

\bibitem[\protect\citeauthoryear{Cao and Bickel}{Cao and
  Bickel}{2022}]{bickel2020}
S. Cao, P.~Bickel,
\newblock Correlations with tailored extremal properties,
\newblock Available at \url{https://arxiv.org/abs/2008.10177v2} (2022).

\bibitem[\protect\citeauthoryear{Chatterjee}{Chatterjee}{2020}]{chatterjee2020}
S. Chatterjee, 
\newblock A new coefficient of correlation,
\newblock {J. Amer. Statist. Ass.  116\/} (2020) 2009--2022.

\bibitem[\protect\citeauthoryear{Cifarelli, Conti, and Regazzini}{Cifarelli
  et~al.}{1996}]{ccr1996}
D. Cifarelli, P.~L. Conti, E.~Regazzini,
\newblock On the asymptotic distribution of a general measure of monotone dependence,
\newblock {Ann. Stat.\/}~{24} (1996) 1386--1399.

\bibitem[\protect\citeauthoryear{Coblenz, Holz, Bauer, Grothe, and
  Koch}{Coblenz et~al.}{2020}]{coblenz2020}
M. Coblenz, S.~Holz, H.~Bauer, O.~Grothe, R.~Koch,
\newblock Modelling fuel injector spray characteristics in jet engines by using vine copulas,
\newblock {J. R. Stat. Soc. Ser. C. Appl. Stat.\/}~{69\/} (2020) 863--886.

\bibitem[\protect\citeauthoryear{Darsow, Nguyen, and Olsen}{Darsow
  et~al.}{1992}]{darsow1992}
W. Darsow, B.~Nguyen, E.~Olsen,
\newblock Copulas and {M}arkov processes,
\newblock {Illinois J. Math.\/}~{36\/} (1992) 600--642.

\bibitem[\protect\citeauthoryear{Deb, Ghosal, and Sen}{Deb
  et~al.}{2020}]{deb2020}
N. Deb, P.~Ghosal, B.~Sen,
\newblock Measuring association on topological spaces using kernels and geometric graphs,
\newblock Available at \url{http://128.84.4.18/abs/2010.01768} (2020).

\bibitem[\protect\citeauthoryear{Dette, Siburg, and Stoimenov}{Dette
  et~al.}{2013}]{siburg2013}
H. Dette, K.~F. Siburg, P.~A. Stoimenov,
\newblock A copula-based non-parametric measure of regression dependence,
\newblock {Scand. J. Statist.\/}~{40\/} (2013), 21--41.

\bibitem[\protect\citeauthoryear{Dette, Van~Hecke, and Volgushev}{Dette
  et~al.}{2014}]{dette2014}
H. Dette, R.~Van~Hecke, S.~Volgushev,
\newblock Some comments on copula-based regression,
\newblock {J. Amer. Statist. Assoc.\/}~{109} (2014) 1319--1324.

\bibitem[\protect\citeauthoryear{Dolati and {\'U}beda-Flores}{Dolati and
  {\'U}beda-Flores}{2006}]{duf2006}
A. Dolati, M.~{\'U}beda-Flores,
\newblock On measures of multivariate concordance,
\newblock {J. Probab. Stat. Sci.\/}~{4\/} (2006) 147--163.

\bibitem[\protect\citeauthoryear{Durante and Fuchs}{Durante and
  Fuchs}{2019}]{sfx2019refl}
F. Durante, S.~Fuchs,
\newblock Reflection invariant copulas,
\newblock {Fuzzy Sets and Systems\/}~{354\/} (2019) 63--73.

\bibitem[\protect\citeauthoryear{Durante and Sempi}{Durante and
  Sempi}{2016}]{fdsempi2016}
F. Durante, C.~Sempi,
\newblock {Principles of Copula Theory},
\newblock CRC Press, Boca Raton FL, 2016.

\bibitem[\protect\citeauthoryear{Fern{\'a}ndez~S{\'a}nchez and
  Trutschnig}{Fern{\'a}ndez~S{\'a}nchez and Trutschnig}{2015}]{fswt2015}
J. Fern{\'a}ndez~S{\'a}nchez, W.~Trutschnig,
\newblock Conditioning-based metrics on the space of multivariate copulas and their interrelation with uniform and levelwise convergence and iterated function systems,
\newblock {J. Theor. Probab.\/}~{28} (2015) 1311--1336.

\bibitem[\protect\citeauthoryear{Fredricks and Nelsen}{Fredricks and
  Nelsen}{2003}]{frednels2003}
G.~A. Fredricks, R.~B. Nelsen,
\newblock The {B}ertino family of copulas,
\newblock In C.~M. Cuadras, J.~Fortiana, J.~A. Rodriguez-Lallena (Eds.),
  {Distributions with Given Marginals and Statistical Modelling}, pp.\ 81--91, Kluwer Academic Publishers, Dordrecht, 2003.

\bibitem[\protect\citeauthoryear{Fuchs}{Fuchs}{2014}]{sfx2014gamma2}
S. Fuchs, 
\newblock Multivariate copulas: Transformations, symmetry, order and measures of concordance,
\newblock {Kybernetika\/}~{50\/} (2014) 725--743.

\bibitem[\protect\citeauthoryear{Fuchs}{Fuchs}{2016}]{sfx2016moc}
S. Fuchs, 
\newblock Copula-induced measures of concordance,
\newblock {Depend. Model.\/}~{4} (2016) 205--214.

\bibitem[\protect\citeauthoryear{Fuchs and McCord}{Fuchs and
  McCord}{2019}]{sfx2019spearman}
S. Fuchs, Y.~McCord,
\newblock On the lower bound of {S}pearman's footrule,
\newblock {Depend. Model.\/}~{7} (2019) 121--129.

\bibitem[\protect\citeauthoryear{Fuchs and Tschimpke}{Fuchs and Tschimpke}{2023}]{sfx2023TP2}
S. Fuchs, M.~Tschimpke,
\newblock Total positivity of copulas from a Markov kernel perspective,
\newblock {J. Math. Anal. Appl. \/}~{518} (2023) Arcticle ID 126629.

\bibitem[\protect\citeauthoryear{Gamboa, Gremaud, Klein, and Lagnoux}{Gamboa
  et~al.}{2020}]{gamboa2020}
F. Gamboa, P.~Gremaud, T.~Klein, A.~Lagnoux,
\newblock Global sensitivity analysis: A novel generation of mighty estimators based on rank statistics,
\newblock {Bernoulli 28\/} (2022) 2345--2374. 
Available at \url{https://arxiv.org/abs/2003.01772}.

\bibitem[\protect\citeauthoryear{Genest, Ne{\v s}lehov{\'a}, and
  Ben~Ghorbal}{Genest et~al.}{2010}]{genebg2010}
C. Genest, J.~Ne{\v s}lehov{\'a}, N.~Ben~Ghorbal,
\newblock Spearman's footrule and {G}ini's gamma: a review with complements,
\newblock {J. Nonparametr. Stat.\/}~{22\/} (2010) 937--954.

\bibitem[\protect\citeauthoryear{Genest, Ne{\v s}lehov{\'a}, and
  Ben~Ghorbal}{Genest et~al.}{2011}]{gene2011}
C. Genest, J.~Ne{\v s}lehov{\'a}, N.~Ben~Ghorbal,
\newblock Estimators based on {K}endall's tau in multivariate copula models,
\newblock {Aust. N. Z. J. Stat.\/}~{53} (2011) 157--177.

\bibitem[\protect\citeauthoryear{Griessenberger, Junker, and
  Trutschnig}{Griessenberger et~al.}{2022}]{fgwt2021}
F. Griessenberger, R.~Junker, W.~Trutschnig,
\newblock On a multivariate copula-based dependence measure and its estimation,
\newblock {Electron. J. Statist.\/}~{16} (2022) 2206--2251.

\bibitem[\protect\citeauthoryear{Han}{Han}{2021}]{bernoulli2021}
F. Han, 
\newblock On extensions of rank correlation coefficients to multivariate spaces,
\newblock {Bernoulli\/}~{28} (2021) 7--11.

\bibitem[\protect\citeauthoryear{Harder and Stadtm{\"u}ller}{Harder and
  Stadtm{\"u}ller}{2014}]{stadtmuller2014}
M. Harder, U.~Stadtm{\"u}ller,
\newblock Maximal non-exchangeability in dimension $d$,
\newblock {J. Multivariate Anal.\/}~{124} (2014) 31--41.

\bibitem[\protect\citeauthoryear{Huang, Deb, and Sen}{Huang
  et~al.}{2020}]{deb2020b}
Z. Huang, N.~Deb, B.~Sen,
\newblock Kernel partial correlation coefficient — a measure of conditional dependence, 
\newblock Available at \url{https://arxiv.org/abs/2012.14804} (2020).

\bibitem[\protect\citeauthoryear{Joe}{Joe}{1990}]{joe1990}
H. Joe,  
\newblock Multivariate concordance,
\newblock {J. Multivariate Anal.\/}~{35} (1990) 12--30.

\bibitem[\protect\citeauthoryear{Joe}{Joe}{1997}]{joe1997}
H. Joe, 
\newblock {Multivariate Models and Dependence Concepts}, 
\newblock Chapman \& Hall, London, 1997.

\bibitem[\protect\citeauthoryear{Junker, Griessenberger, and Trutschnig}{Junker
  et~al.}{2020}]{fgwt2020}
R. Junker, F.~Griessenberger, W.~Trutschnig,
\newblock Estimating scale-invariant directed dependence of bivariate distributions, 
\newblock {Comput. Statist. Data Anal.\/}~{153} (2020) Article ID 107058.

\bibitem[\protect\citeauthoryear{Kallenberg}{Kallenberg}{2002}]{kallenberg2002}
O. Kallenberg, 
\newblock {Foundations of {M}odern Probability},
\newblock Springer, New York, 2002.

\bibitem[\protect\citeauthoryear{Kamnitui, Fern{\'a}ndez-S{\'a}nchez, and
  Trutschnig}{Kamnitui et~al.}{2018}]{kamnitui2018}
N. Kamnitui, J.~Fern{\'a}ndez-S{\'a}nchez, W.~Trutschnig,
\newblock Maximum asymmetry of copulas revisited,
\newblock {Depend. Model.\/}~{6} (2018) 47--62.

\bibitem[\protect\citeauthoryear{Karger, Conrad, B{\"o}hner, Kawohl, Kreft,
  Soria-Auza, Zimmermann, Linder, and Kessler}{Karger
  et~al.}{2017}]{karger2017}
D. Karger, O.~Conrad, J.~B{\"o}hner, T.~Kawohl, H.~Kreft, R.~Soria-Auza, N.~Zimmermann, H.~Linder, M.~Kessler,
\newblock Climatologies at high resolution for the earth's land surface areas,
\newblock {Sci. Data\/}~{4} (2017) Article ID 170122.

\bibitem[\protect\citeauthoryear{Karger, Conrad, B{\"o}hner, Kawohl, Kreft,
  Soria-Auza, Zimmermann, Linder, and Kessler}{Karger
  et~al.}{2018}]{karger2018}
D. Karger, O.~Conrad, J.~B{\"o}hner, T.~Kawohl, H.~Kreft, R.~Soria-Auza, N.~Zimmermann, H.~Linder, M.~Kessler,
\newblock Data from: Climatologies at high resolution for the earth's land surface areas (2018) [dataset].

\bibitem[\protect\citeauthoryear{Kasper, Fuchs, and Trutschnig}{Kasper
  et~al.}{2021}]{sfx2021weak}
T. Kasper, S.~Fuchs, W.~Trutschnig,
\newblock On weak conditional convergence of bivariate {A}rchimedean and {E}xtreme {V}alue copulas, and consequences to nonparametric estimation,
\newblock {Bernoulli\/}~{27} (2021) 2217--2240.

\bibitem[\protect\citeauthoryear{Klenke}{Klenke}{2008}]{klenke2008}
A. Klenke,  
\newblock {Wahrscheinlichkeitstheorie},
\newblock Springer, Heidelberg, 2008.

\bibitem[\protect\citeauthoryear{Lancaster}{Lancaster}{1963}]{lancaster1963}
H.~O. Lancaster,
\newblock Correlation and complete dependence of random variables,
\newblock {Ann. Math. Statist.\/}~{34\/} (1963) 1315--1321.

\bibitem[\protect\citeauthoryear{Li, Scarsini, and Shaked}{Li
  et~al.}{1996}]{lss1996}
X. Li, M.~Scarsini, M.~Shaked,
\newblock Linkages: A tool for the construction of multivariate distributions with given nonoverlapping multivariate marginals,
\newblock {J. Multivariate Anal.\/}~{56} (1996) 20--41.

\bibitem[\protect\citeauthoryear{Lopez-Paz, Hennig, and
  Sch{\"o}lkopf}{Lopez-Paz et~al.}{2013}]{lopez2013}
D. Lopez-Paz, P.~Hennig, B.~Sch{\"o}lkopf,
\newblock The randomized dependence coefficient,
\newblock {Advances in Neural Information Processing Systems\/} (2013) 1--9.

\bibitem[\protect\citeauthoryear{McDiarmid}{McDiarmid}{1989}]{mcdiarmid1989}
C. McDiarmid, 
\newblock On the method of bounded differences,
\newblock In J.~Siemons (Ed.), {Surveys in Combinatorics}, pp.\  144--188,
  Cambridge University Press, 1989.

\bibitem[\protect\citeauthoryear{Mroz, Fern{\'a}ndez-S{\'a}nchez, Fuchs, and
  Trutschnig}{Mroz et~al.}{2021}]{sfx2021endo}
T. Mroz, J.~Fern{\'a}ndez-S{\'a}nchez, S.~Fuchs, W.~Trutschnig,
\newblock On distributions with fixed marginals maximizing the joint or the prior default probability, estimation, and related results,
\newblock {J. Statist. Plann. Inference 223} (2023) 33--52.

\bibitem[\protect\citeauthoryear{Mroz, Fuchs, and Trutschnig}{Mroz
  et~al.}{2021}]{sfx2021vine}
T. Mroz, S.~Fuchs, W.~Trutschnig, 
\newblock How simplifying and flexible is the simplifying assumption in pair-copula constructions – analytic answers in dimension three and a glimpse beyond,
\newblock {Electron. J. Statist.\/}~{15\/} (2021) 1951--1992.

\bibitem[\protect\citeauthoryear{Nelsen}{Nelsen}{2006}]{nelsen2006}
R.~B. Nelsen,
\newblock {An Introduction to Copulas},
\newblock Springer, New York, 2006.

\bibitem[\protect\citeauthoryear{Nelsen}{Nelsen}{2007}]{nel2007}
R.~B. Nelsen, 
\newblock Extremes of nonexchangeability,
\newblock {Statist. Papers\/}~{48} (2007) 329--336.

\bibitem[\protect\citeauthoryear{Noh, El~Ghouch, and Bouzmarni}{Noh
  et~al.}{2013}]{noh2013}
H. Noh, A.~El~Ghouch, T.~Bouzmarni,
\newblock Copula-based regression estimation and inference,
\newblock {J. Amer. Statist. Assoc.\/}~{108} (2013) 676--688.

\bibitem[\protect\citeauthoryear{Reshef, Reshef, Finucane, Grossman, McVean,
  Turnbaugh, Lander, Mitzenmacher, and Sabeti}{Reshef
  et~al.}{2011}]{reshef2011}
D. Reshef, Y.~Reshef, H.~Finucane, S.~Grossman, G.~McVean, P.~Turnbaugh, E.~Lander, M.~Mitzenmacher, P.~Sabeti,
\newblock Detecting novel associations in large data sets,
\newblock {Science\/}~{334\/} (2011) 1518--1524.

\bibitem[\protect\citeauthoryear{Schmid and Schmidt}{Schmid and
  Schmidt}{2007}]{schmidt2006}
F. Schmid, R.~Schmidt,
\newblock Multivariate extensions of {S}pearman’s rho and related statistics,
\newblock {Stat. Probab. Lett.\/}~{77} (2007) 407--416.

\bibitem[\protect\citeauthoryear{Shi, Drton, and Han}{Shi
  et~al.}{2021}]{shi2021normal}
H. Shi, M.~Drton, F.~Han,
\newblock On {A}zadkia-{C}hatterjee's conditional dependence coefficient,
\newblock Available at \url{https://arxiv.org/abs/2108.06827v1} (2021).

\bibitem[\protect\citeauthoryear{Shih and Emura}{Shih and
  Emura}{2021}]{emura2021}
J.-H. Shih, T.~Emura,
\newblock On the copula correlation ratio and its generalization,
\newblock {J. Multivariate Anal.\/}~{182} (2021) Article ID 104708.

\bibitem[\protect\citeauthoryear{Strothmann, Dette, and Siburg}{Strothmann
  et~al.}{2022}]{strothmann2022}
C. Strothmann, H.~Dette, K.~Siburg,
\newblock Rearranged dependence measures,
\newblock Available at \url{https://arxiv.org/abs/2201.03329v1} (2022).

\bibitem[\protect\citeauthoryear{Sungur}{Sungur}{2005}]{sungur2005}
E. Sungur, 
\newblock A note on directional dependence in regression setting,
\newblock {Comm. Statist. Theory Methods\/}~{34} (2005) 1957--1965.

\bibitem[\protect\citeauthoryear{Sz{\'e}kely, Rizzo, and Bakirov}{Sz{\'e}kely
  et~al.}{2007}]{szekely2007}
G. Sz{\'e}kely, M.~Rizzo, N.~Bakirov,
\newblock Measuring and testing dependence by correlation of distances,
\newblock {Ann. Statist.\/}~{35} (2007) 2769--2794.

\bibitem[\protect\citeauthoryear{Taylor}{Taylor}{2016}]{tay2016}
M.~D. Taylor, 
\newblock Multivariate measures of concordance for copulas and their marginals,
\newblock {Depend. Model.\/}~{4} (2016) 224--236.

\bibitem[\protect\citeauthoryear{Trutschnig}{Trutschnig}{2011}]{wt2011}
W. Trutschnig,
\newblock On a strong metric on the space of copulas and its induced dependence measure,
\newblock {J. Math. Anal. Appl.\/}~{384\/} (2011) 690--705.

\bibitem[\protect\citeauthoryear{Trutschnig and
  Fern{\'a}ndez-S{\'a}nchez}{Trutschnig and
  Fern{\'a}ndez-S{\'a}nchez}{2012}]{wt2012}
W. Trutschnig, J.~Fern{\'a}ndez-S{\'a}nchez,
\newblock Idempotent and multivariate copulas with fractal support,
\newblock {J. Statist. Plann. Inference\/}~{142} (2012) 3086--3096.
\end{thebibliography}
%\bibliographystyle{chicago}

%%%%%%%%%%%%%%%%%%%%%%%%%%%%%%%%%%%%%%%%%%%%%%%%%%%%%%%%%%%%%%%%%%%%%%%%%%%%%%%%%%%%%%%%%%%%%%%%%%%%%%%%%%%%%%%%%%
%%%%%%%%%%%%%%%%%%%%%%%%%%%%%%%%%%%%%%%%%%%%%%%%%%%%%%%%%%%%%%%%%%%%%%%%%%%%%%%%%%%%%%%%%%%%%%%%%%%%%%%%%%%%%%%%%%
%%%%%%%%%%%%%%%%%%%%%%%%%%%%%%%%%%%%%%%%%%%%%%%%%%%%%%%%%%%%%%%%%%%%%%%%%%%%%%%%%%%%%%%%%%%%%%%%%%%%%%%%%%%%%%%%%%

\appendix

\section{Proofs} 
%\ref{Sect.phi}

\noindent
{\bf Proof of Theorem \ref{phi.CompleteDep}.} 
We first prove (i). 
Therefore, assume that $Y$ is completely dependent on $\XXX$.
Then there exists some $\mu_{A^{1:d}}$-$\lambda$-preserving map $h$ such that 
$ K_A (\uuu, [0,v])
	= \mathds{1}_{[0,v]} (h(\uuu)) $ for $\mu_{A^{1:d}}$-almost all $\uuu \in \I^d$ and every $v \in \I$ (see \cite{fgwt2021}).
In this case
$ K_A (\uuu, [0,s]) \, K_A (\uuu, [0,t])
	= \mathds{1}_{[0,s]} (h(\uuu)) \mathds{1}_{[0,t]} (h(\uuu))
	= \mathds{1}_{[0,M(s,t)]} (h(\uuu))
	= K_A (\uuu, [0,M(s,t)]) $
 and hence 
$$
  \psi(A)(s,t)
	%= \int_{[0,1]^d} K_A (\uuu, [0,s]) \, K_A (\uuu, [0,t]) \; \mathrm{d} \mu_{A^d} (\uuu)
	= \int_{\I^d} K_A (\uuu, [0,M(s,t)]) \; \mathrm{d} \mu_{A^{1:d}} (\uuu)
	= M(s,t),
	\qquad (s,t) \in \I^2. 
$$
Now, assume that $\psi(A) = M$.
Then
$ v
	  =  M(v,v)
	  =  \psi(A)(v,v)
	  =  \int_{\I^d} K_A (\uuu, [0,v])^2 \; \mathrm{d} \mu_{A^{1:d}} (\uuu)
	\leq \int_{\I^d} K_A (\uuu, [0,v]) \; \mathrm{d} \mu_{A^{1:d}} (\uuu) 
	  =  v $
for all $v \in \I$ which implies
$ \int_{\I^d} K_A (\uuu, [0,v]) - K_A (\uuu, [0,v])^2 \; \mathrm{d} \mu_{A^{1:d}} (\uuu) = 0 $.
Since the integrand is positive we obtain
$K_A (\uuu, [0,v]) - K_A (\uuu, [0,v])^2 = 0$
or, equivalently, either $K_A (\uuu, [0,v]) = 0$ or $K_A (\uuu, [0,v]) = 1$ for $\mu_{A^{1:d}}$-almost every $\uuu \in \I^d$.
Since $v \mapsto K_A (\uuu, [0,v])$ is increasing, 
it can be shown in exactly the same manner as in \cite[Lemma 12]{wt2011} that there exists some $\mu_{A^{1:d}}-\lambda$ preserving transformation $h: \I^d \to \I$ such that 
$K_A (\uuu, F) 
%= \mathds{1}_{[h(\uuu),1]} (v) 
= \mathds{1}_{F} (h(\uuu))$ is a regular conditional distribution of $A$. 
This proves (i). 
\\
We now prove (ii).
Therefore, assume that $Y$ and $\XXX$ are independent.
Then 
$ \psi(A)(s,t) 
	%= \int_{\I^d} K_A (\uuu, [0,s]) \, K_A (\uuu, [0,t]) \; \mathrm{d} \mu_{A^{1:d}} (\uuu)
	= \int_{\I^d} st \; \mathrm{d} \mu_{A^{1:d}} (\uuu)
	= \Pi(s,t) $
for all $(s,t) \in \I^2$ and hence (b), which then implies (c). 
Now, assume that $\delta_{\psi(A)} = \delta_{\Pi}$.
Applying H{\"o}lder's inequality yields
$$
  v^2 
	  =  \delta_{\Pi}(v)
	  =  \delta_{\psi(A)}(v)
	  =  \int_{\I^d} K_A (\uuu, [0,v])^2 \; \mathrm{d} \mu_{A^{1:d}} (\uuu)
	\geq \left( \int_{\I^d} K_A (\uuu, [0,v]) \; \mathrm{d} \mu_{A^{1:d}} (\uuu) \right)^2
	  =  v^2.
$$
Hence we have
$ 0
	= \int_{\I^d} K_A (\uuu, [0,v])^2 - v^2 \; \mathrm{d} \mu_{A^{1:d}} (\uuu)
	= \int_{\I^d} \big( K_A (\uuu, [0,v]) - v \big)^2 \; \mathrm{d} \mu_{A^{1:d}} (\uuu) $
for all $v \in \I$.
Since the integrand is positive we obtain
$K_A (\uuu, [0,v]) = v$
for $\mu_{A^{1:d}}$-almost every $\uuu \in \I^d$ and all $v \in \I$, hence $A(\uuu,v) = A^{1:d}(\uuu) \, v$ for all $(\uuu,v) \in \I^{d}\times\I$.
This proves (ii).
\qedsymbol

\bigskip\noindent
{\bf Proof of Corollary \ref{phi.Diagonal}.} 
Applying H{\"o}lder's inequality directly shows that the diagonal of $\psi(A)$ exceeds the diagonal of $\Pi$, 
i.e., $\delta_\Pi(t) := t^2 = \Pi(t,t) \leq \psi(A)(t,t)$ for all $t \in \I$.
According to \cite[Proposition 1.2]{frednels2003}, the Bertino copula 
$ B_\delta(u,v)
	:= \min\{u,v\}-\min\{t - \delta(t) \, : \, t \in [u,v]\} $ 
has diagonal $\delta$ and serves as a lower bound within the class of copulas having diagonal $\delta$.
Hence, the Bertino copula $B_{\delta_\Pi}$ with
$$
  B_{\delta_\Pi} (s,t) 
	= \min(s,t) - \min\{ p - p^2 \, : \, p \in [s,t]\}
	= \begin{cases}
	    \min(s,t)^2, 
		  & s+t \leq 1, 
		  \\
		  \min(s,t) - \max(s,t) + \max(s,t)^2,
		  & s+t > 1,
	  \end{cases}
$$
(see \cite[Example 1.5]{frednels2003})
is a lower bound within the class of copula having diagonal $\delta_\Pi$.
Since $B_{\delta} \leq B_{\delta^\prime}$ whenever $\delta \leq \delta^\prime$ we may conclude that
$B_{\delta_\Pi} \leq \psi(A)$ for all $A \in \mathcal{C}^{d+1}$.
This proves the assertion.
\qedsymbol

\bigskip\noindent
{\bf Proof of Theorem \ref{phi.InformationGain}.}
Due to Corollary \ref{phi.Trans}, 
we may assume that $L = \{1,\dots,l\}$ for some $l \in \{1,\dots,d-1\}$.
Then, with the notation used at the beginning of Section \ref{Sect.phi}
%we consider the $(d+2)$-dimensional random vector $(\UUU,V,V^{\prime})$ with uniform univariate marginals such that $(\UUU,V) \sim A$ and %$(\UUU,V^{\prime}) \sim A$ 
%and assume that $V$ and $V^{\prime}$ are conditionally independent conditioned on $\UUU$.
%Then, 
\begin{eqnarray*}
  \psi(A^{1:l,d+1}) (t,t)
	& = & E \big( E (\mathds{1}_{[0,t]} \circ V \, | \, (U_1,\dots,U_{l}))^2 \big),
	\\
	\psi(A^{1:(l+1),d+1}) (t,t)
	& = & E \big( E (\mathds{1}_{[0,t]} \circ V \, | \, (U_1,\dots,U_{l},U_{l+1}))^2 \big),
\end{eqnarray*}
for all $t \in \I$.
Considering the right hand side as conditional expectation 
\begin{eqnarray*}
  \psi(A^{1:l,d+1}) (t,t)
	& = & E ((\mathds{1}_{[0,t]} \circ V)^2)  - E \Big( \Big( \mathds{1}_{[0,t]} \circ V - E \big( \mathds{1}_{[0,t]} \circ V \, | \, (U_1,\dots,U_{l}) \big) \Big)^2 \Big),
	\\*
	\psi(A^{1:(l+1),d+1}) (t,t)
	& = &  E ((\mathds{1}_{[0,t]} \circ V)^2)  - E \Big( \Big( \mathds{1}_{[0,t]} \circ V - E \big( \mathds{1}_{[0,t]} \circ V \, | \, (U_1,\dots,U_{l},U_{l+1}) \big) \Big)^2 \Big),
\end{eqnarray*}
and applying Hilbert's projection theorem yields
$$
  E \Big( \Big( \mathds{1}_{[0,t]} \circ V - E \big( \mathds{1}_{[0,t]} \circ V \, | \, (U_1,\dots,U_{l},U_{l+1}) \big) \Big)^2 \Big)
	\leq E \Big( \Big( \mathds{1}_{[0,t]} \circ V - E \big( \mathds{1}_{[0,t]} \circ V \, | \, (U_1,\dots,U_{l}) \big) \Big)^2 \Big).
$$
Hence 
$ \psi(A^{1:l,d+1}) (t,t)
	\leq \psi(A^{1:(l+1),d+1}) (t,t) $
for all $t \in \I$.
This proves the assertion.
\qedsymbol

\bigskip
\begin{lemma}{} \label{Lemma.CI}
For a $(d+1)$-dimensional copula $A$ and $k \in \{1,\dots,d-1\}$, the following statements are equivalent:
\begin{enumerate}[(i)]
\item
The identity
$ K_{A} (\uuu,[{\bf 0},\vvv] \times [0,t])
	= K_{A^{1:k,d+1}} (\uuu, [0,t]) \, K_{A^{1:d}} (\uuu, [{\bf 0},\vvv]) $
holds for all $(\vvv,t) \in \I^{d-k} \times \I$ and $\mu_{A^{1:k}}$-a.e. $\uuu \in \I^k$.
\item
The identity
$ K_{A} ((\uuu,\vvv),[0,t]) = K_{A^{1:k,d+1}} (\uuu, [0,t]) $ 
holds for all $t \in \I$ and $\mu_{A^{1:d}}$-a.e. $(\uuu,\vvv) \in \I^{k} \times \I^{d-k}$.
\item
$\psi(A) (t,t) = \psi(A^{1:k,d+1}) (t,t)$ for all $t \in \I$.
\end{enumerate}
In addition, (a) even implies $\psi(A) = \psi(A^{1:k,d+1})$.
\end{lemma}{}
\noindent
{\bf Proof.}
Fix $k \in \{1,\dots,d-1\}$.
The equivalence of (i) and (ii) follows from the identities
$$
  \int_{[{\bf 0},\uuu]} K_{A} (\ppp,[{\bf 0},\vvv] \times [0,t]) \; \mathrm{d} \mu_{A^{1:k}}(\ppp)
	= A (\uuu,\vvv, t)
	= \int_{[{\bf 0},\uuu] \times [{\bf 0},\vvv]} K_{A} ((\ppp,\qqq),[0,t]) \; \mathrm{d} \mu_{A^{1:d}} (\ppp,\qqq),
	\qquad (\uuu,\vvv,t) \in \I^{k} \times \I^{d-k} \times \I,
$$
and
\begin{eqnarray*}
	\int_{[{\bf 0},\uuu] \times [{\bf 0},\vvv]} K_{A^{1:k,d+1}} (\ppp, [0,t]) \; \mathrm{d} \mu_{A^{1:d}} (\ppp,\qqq)
	& = & \int_{[{\bf 0},\uuu] \times [{\bf 0},\vvv]} K_{A^{1:k,d+1}} (\ppp, [0,t]) \, K_{A^{1:d}} (\ppp, \mathrm{d}\qqq) \, \mathrm{d} \mu_{A^{1:k}}(\ppp)
	\\
	& = &  \int_{[{\bf 0},\uuu]} K_{A^{1:k,d+1}} (\ppp, [0,t]) \, K_{A^{1:d}} (\ppp, [{\bf 0},\vvv]) \; \mathrm{d} \mu_{A^{1:k}}(\ppp),
	\qquad (\uuu,\vvv,t) \in \I^{k} \times \I^{d-k} \times \I,
\end{eqnarray*}
by applying Radon-Nikodym theorem.
Now, assume that (ii) holds.
Then, for all $(s,t) \in \I^2$, 
\begin{align*}
  \psi(A) (s,t)
	& = \int_{\I^d} K_{A} ((\uuu,\vvv), [0,s]) \, K_{A} ((\uuu,\vvv), [0,t]) \; \mathrm{d} \mu_{A^{1:d}} (\uuu,\vvv)
	  = \int_{\I^d} K_{A^{1:k,d+1}} (\uuu, [0,s]) \, K_{A^{1:k,d+1}} (\uuu, [0,t]) \; \mathrm{d} \mu_{A^{1:d}} (\uuu,\vvv)
	\\
	& =  \int_{\I^k} K_{A^{1:k,d+1}} (\uuu, [0,s]) \, K_{A^{1:k,d+1}} (\uuu, [0,t]) 
				\int_{\I^{d-k}} K_{A^{1:d}} (\uuu, \mathrm{d} \vvv) \; \mathrm{d} \mu_{A^{1:k}} (\uuu)
				\\
	& =  \int_{\I^k} K_{A^{1:k,d+1}} (\uuu, [0,s]) \, K_{A^{1:k,d+1}} (\uuu, [0,t]) \, K_{A^{1:d}} (\uuu, \I^{d-k}) \; \mathrm{d} \mu_{A^{1:k}} (\uuu)
	\\
	&	=  \int_{\I^k} K_{A^{1:k,d+1}} (\uuu, [0,s]) \, K_{A^{1:k,d+1}} (\uuu, [0,t]) \; \mathrm{d} \mu_{A^{1:k}} (\uuu)
	  =  \psi(A^{1:k,d+1}) (s,t).
\end{align*}
This proves (iii) and the additional result. 
Finally, assume that (iii) holds.
Then, we have 
$ \psi(A^{1:k,d+1}) (t,t)
	= \psi(A) (t,t) $
and hence, for all $t \in \I$, 
\begin{eqnarray*}
  0
	%& = & \psi(A) (t,t) - \psi(A^{1,d+1}) (t,t)
	%\\
	& = & \int_{\I^{d}} K_{A} ((\uuu,\vvv), [0,t])^2 - K_{A^{1:k,d+1}} (\uuu, [0,t])^2 \; \mathrm{d} \mu_{A^{1:d}} (\uuu,\vvv).
\end{eqnarray*}
Since, for every $t \in \I$,
$ K_{A^{1:k,d+1}} (\uuu, [0,t])
	= \int_{\I^{d-k}} K_{A} ((\uuu,\vvv), [0,t]) \; K_{A^{1:d}}(\uuu, \mathrm{d} \vvv) $
for $\mu_{A^{1:k}}$-a.e. $\uuu \in \I^k$, and hence
\begin{eqnarray*}
   \lefteqn{\int_{\I^{d}} K_{A^{1:k,d+1}} (\uuu, [0,t]) \, K_{A} ((\uuu,\vvv), [0,t]) \; \mathrm{d} \mu_{A^{1:d}} (\uuu,\vvv)
	 = \int_{\I^k} \int_{\I^{d-k}} K_{A^{1:k,d+1}} (\uuu, [0,t]) \, K_{A} ((\uuu,\vvv), [0,t]) \; K_{A^{1:d}}(\uuu, \mathrm{d} \vvv) \; \mathrm{d} \mu_{A^{1:k}}(\uuu)}
	\\
	& = & \int_{\I^k} K_{A^{1:k,d+1}} (\uuu, [0,t]) \int_{\I^{d-k}} K_{A} ((\uuu,\vvv), [0,t]) \; K_{A^{1:d}}(\uuu, \mathrm{d} \vvv) \; \mathrm{d} \mu_{A^{1:k}}(\uuu) 
	 = \int_{\I^k} K_{A^{1:k,d+1}} (\uuu, [0,t])^2  \; \mathrm{d} \mu_{A^{1:k}}(\uuu), 
\end{eqnarray*}
we finally obtain
$ 0
	= \int_{\I^{d}} \Big( K_{A} ((\uuu,\vvv), [0,t]) - K_{A^{1:k,d+1}} (\uuu, [0,t]) \Big)^2 \; \mathrm{d} \mu_{A^{1:d}} (\uuu,\vvv) $
and hence $K_{A} ((\uuu,\vvv), [0,t]) = K_{A^{1:k,d+1}} (\uuu, [0,t])$ for all $t \in \I$ and $\mu_{A^{1:d}}$-a.e. $(\uuu,\vvv) \in \I^{k} \times \I^{d-k}$,
which proves (ii).
\qedsymbol

\bigskip\noindent
{\bf Proof of Corollary \ref{phi.Trans}.}  
Since $\zeta$ is bijective, 
the Markov kernel of $A_{(\zeta,{\rm id})}$ satisfies 
$ K_{A_{(\zeta,{\rm id})}}(\uuu,[0,v])
	= K_{A}(\zeta^{-1}(\uuu),[0,v]) $
for $(\mu_{A^{1:d}})^{\zeta}$-almost all $\uuu \in \I^{d}$ and all $v \in \I$ and hence 
\begin{align*}
  \psi(A_{(\zeta,{\rm id})}) (s,t)
	& = \int_{\I^d} K_{A_{(\zeta,{\rm id})}} (\uuu, [0,s]) \, K_{A_{(\zeta,{\rm id})}} (\uuu, [0,t]) \; \mathrm{d} \mu_{(A_{(\zeta,{\rm id})})^{1:d}} (\uuu)
	  = \int_{\I^d} K_{A} (\zeta^{-1}(\uuu), [0,s]) \, K_{A} (\zeta^{-1}(\uuu), [0,t]) \; \mathrm{d} (\mu_{A^{1:d}})^{\zeta} (\uuu)
	\\
	& = \int_{\I^d} K_{A} (\zeta^{-1}(\zeta(\uuu)), [0,s]) \, K_{A} (\zeta^{-1}(\zeta(\uuu)), [0,t]) \; \mathrm{d} \mu_{A^{1:d}} (\uuu)
	  %= \int_{\I^d} K_{A} (\uuu, [0,s]) \, K_{A} (\uuu, [0,t]) \; \mathrm{d} \mu_{A^{1:d}} (\uuu)
	  = \psi(A)(s,t),
\end{align*}
which proves the result.
\qedsymbol

\bigskip\noindent
{\bf Proof of Theorem \ref{zeta.Representation}.} 
Denote by $H$ the distribution function of $(\XXX,Y)$.
Since $(\XXX,Y)$ has continuous marginals, the random vector $({\bf F} (\XXX), G(Y))$ has distribution function $A$.
%and $P^\XXX = (P^{{\bf F}(\XXX)})^{{\bf F}^\leftarrow}$.
By change of measure we first obtain
$ \int_{\mathbb{R}} {\rm var} (\mathds{1}_{\{Y \geq y\}}) \; \mathrm{d} P^{Y}(y)
	= \int_{\I} P(Y \geq G^\leftarrow (v)) [1 - P(Y \geq G^\leftarrow (v))] \; \mathrm{d} \lambda(v)
	= \int_{\I} (1-v)\, v \; \mathrm{d} \lambda(v)
	= \frac{1}{6} $
and \cite[Lemma 1]{sfx2021endo} gives
\begin{align*}
  \int_{\mathbb{R}} {\rm var} (P(Y \geq y \, | \, \XXX)) \; \mathrm{d} P^{Y}(y)
	& = \int_{\I} E\big( P(Y \geq G^\leftarrow (v) \, | \, \XXX)^2 \big) 
								- P(Y \geq  G^\leftarrow (v))^2 \; \mathrm{d} \lambda(v)
	\\*
	& = \int_{\I} E\big( K_H (\XXX, [G^\leftarrow (v),\infty) )^2 \big) 
								- (1-v)^2 \; \mathrm{d} \lambda(v)
	  = \int_{\I} E\big( K_A ({\bf F}(\XXX), [v,1] )^2 \big) 
								- (1-v)^2 \; \mathrm{d} \lambda(v)
	\\
	& = \int_{\I} \int_{\I^d} K_A (\uuu, [v,1] )^2 \; \mathrm{d} \mu_{A^{1:d}} (\uuu) \, \mathrm{d} \lambda(v) - \frac{1}{3}
	  = \int_{\I} \int_{\I^d} K_A (\uuu, [0,v))^2 \; \mathrm{d} \mu_{A^{1:d}} (\uuu) \, \mathrm{d} \lambda(v) - \frac{1}{3}
	\\
	& = \int_{\I} \int_{\I^d} K_A (\uuu, [0,v])^2 \; \mathrm{d} \mu_{A^{1:d}} (\uuu) \, \mathrm{d} \lambda(v) - \frac{1}{3}
	  = \int_{\I} \psi(A)(t,t) \; \mathrm{d} \lambda(t) - \frac{1}{3}.
\end{align*}
Finally, 
$ T(Y,\XXX)
	= 6 \, \int_{\I} \psi(A)(t,t) \; \mathrm{d} \lambda(t) - 2
	= \phi(\psi(A)) $
which proves the assertion.
\qedsymbol

\bigskip\noindent
{\bf Proof of Theorem \ref{R2.Representation.II}.} 
Since
$ E (E (Y \, | \, \XXX)^2)
	= E \big( E (Y \, | \, \XXX) \, E (Y^\prime \, | \, \XXX)\big)
	= E (E (Y \, Y^\prime \, | \, \XXX))
	= E (Y \, Y^\prime) $
we obtain
\begin{align*}
  R^2 (Y,\XXX)
	& = \frac{{\rm var} ( E (Y \, | \, \XXX) )}{{\rm var} (Y)}
	= \frac{E (E (Y \, | \, \XXX)^2) - E(Y)^2}{{\rm var} (Y)}   
	= \frac{E (Y \, Y^\prime)- E(Y)E(Y^\prime)}{\sqrt{{\rm var} (Y)}\sqrt{{\rm var} (Y^\prime)}}
	= \varrho_P (Y,Y^\prime).
\end{align*}
This proves the result.
\qedsymbol

\bigskip\noindent
{\bf Proof of Corollary \ref{R2.Properties}.} 
It remains to prove property (ii) which is immediate from the identity
$R^2(V,\UUU)$
$= 12 \int_{\I^d} \big( E(V | \UUU=\uuu) - \frac{1}{2} \big)^2 \; \mathrm{d} \mu_{A^{1:d}} (\uuu) .
	%= 12 \int_{\I^d} \left( \int_{\I} K_A(\uuu,[0,v]) \; \mathrm{d} \lambda (v) - \frac{1}{2} \right)^2 \; \mathrm{d} \mu_{A^{1:d}} (\uuu)
	$
This proves the assertion.
\qedsymbol

\bigskip\noindent
{\bf Proof of Theorem \ref{Q.Properties}.} 
It remains to prove the identity which is immediate from
\begin{align*}
  \lefteqn{2 \, \int_{\I} \int_{\I^d} \big( K_A(\uuu,[0,t]) - (1 - K_A(\uuu,[0,1-t])) \big)^2 \, \mathrm{d} \mu_{A^{1:d}} (\uuu) \mathrm{d} \lambda(t)}
	\\
	& = 2 \, \int_{\I} \int_{\I^d} K_A(\uuu,[0,t])^2 - 2 \, K_A(\uuu,[0,t]) + 2 \, K_A(\uuu,[0,t]) K_A(\uuu,[0,1-t]) 
	\\
	&   \qquad\qquad + 1 - 2\,K_A(\uuu,[0,1-t]) + K_A(\uuu,[0,1-t])^2 \, \mathrm{d} \mu_{A^{1:d}} (\uuu) \mathrm{d} \lambda(t)
	\\
	& = 2 \, \int_{\I} \int_{\I^d} K_A(\uuu,[0,t])^2 + 2 \, K_A(\uuu,[0,t]) K_A(\uuu,[0,1-t]) + K_A(\uuu,[0,1-t])^2 
				\; \mathrm{d} \mu_{A^{1:d}} (\uuu) \mathrm{d} \lambda(t) - 2
	\\
	& = 2 \, \int_{\I} \psi(A)(t,t) + 2 \psi(A)(t,1-t) + \psi(A)(1-t,1-t) \; \mathrm{d} \lambda(t) - 2
	\\
	& = 4 \, \int_{\I} \psi(A)(t,t) + \psi(A)(t,1-t) \; \mathrm{d} \lambda(t) - 2
	  = \gamma (\psi(A))
	  = Q(Y,\XXX).
\end{align*}
This proves the result.
\qedsymbol

%\label{App.Est.}
\bigskip
For proving consistency of \eqref{Estimation.Estimate.Alternative} (Theorem \ref{Consistency.Thm2}) we use a modification of $D_n$ for which we choose the usual normalization of the ranks.
For $(s,t) \in \I^2$, define 
\begin{equation} \label{Estimation.Estimate}
  C_n (s,t) 
	:= 
	%\frac{1}{n} \, \sum_{k=1}^{n} \mathds{1}_{[0,s]} (R_k/n) \mathds{1}_{[0,t]} (R_{N(k)}/n)		
	%\\
	%&  = & 
	\frac{1}{n} \, \sum_{k=1}^{n} \mathds{1}_{[0,s]} (G_n(Y_k)) \mathds{1}_{[0,t]} (G_n(Y_{N(k)})),	
\end{equation}
where $G_n$ denotes the empirical distribution function of $Y_1,\dots,Y_n$, i.e.,
$ G_n(y)
	= \frac{1}{n} \; \sum_{k=1}^{n} \mathds{1}_{(-\infty,y]} (Y_k) $.
Due to Lemma \ref{Lemma.Chatterjee.1} below, the estimators $C_n$ and $D_n$ are asymptotically equivalent 
(see also (\ref{Est.Asymtotic.Eq}) below).

For a realization $\XXX_1, \dots,\XXX_n$ and each $i \in \{1,\dots,n\}$, let 
$K_{n,i}$ be the number of $j$ such that $\XXX_i$ is the nearest neighbour of $\XXX_j$.
%(not necessarily the randomly chosen one) and $\XXX_j \neq \XXX_i$.
The following result is given in \cite{chatterjee2021}:

\begin{lemma}{} \label{Lemma.Chatterjee.1} \cite[Lemma 11.4]{chatterjee2021} \\
There exists a constant $c(d)$ such that 
$ K_{n,1} \leq c(d) $.
\end{lemma}{}
\noindent
Notice that the upper bound $c(d)$ used in Lemma \ref{Lemma.Chatterjee.1} only depends on dimension $d$.

The following two results are key for proving consistency of \eqref{Estimation.Estimate}:

\begin{lemma}{} \label{thm.convergence.1}
For every $(s,t) \in \I^2$ we have
$\lim_{n \to \infty} E(C_n (s,t)) = \psi(A) (s,t)$.
\end{lemma}{}
\noindent
{\bf Proof.}
For $(s,t) \in \I^2$, we define 
$$
  C_n^\ast (s,t) 
	:= \frac{1}{n} \, \sum_{k=1}^{n} \mathds{1}_{[0,s]} (G(Y_k)) \mathds{1}_{[0,t]} (G(Y_{N(k)})),
$$
and show that 
\begin{enumerate}
\item $\lim_{n \to \infty} E (|C_n(s,t) - C_n^\ast(s,t)|) = 0$, and 
\item $\lim_{n \to \infty} E(C_n^\ast(s,t)) = \psi(A) (s,t)$.
\end{enumerate}
We first prove 1.
For $(s,t) \in \I^2$, we have
\begin{eqnarray*}
  | C_n (s,t)  - C_n^\ast (s,t) |
	& \leq & \frac{1}{n} \, \sum_{k=1}^{n} 
						\big| \mathds{1}_{[0,s]} (G_n(Y_k)) \mathds{1}_{[0,t]} (G_n(Y_{N(k)})) 
									- \mathds{1}_{[0,s]} (G(Y_k)) \mathds{1}_{[0,t]} (G(Y_{N(k)})) \big|
	\\
	& \leq & \frac{1}{n} \, \sum_{k=1}^{n} 
						\big| \mathds{1}_{[0,s]} (G_n(Y_k)) - \mathds{1}_{[0,s]} (G(Y_k)) \big|
						+ \frac{1}{n} \, \sum_{k=1}^{n}  
						\big| \mathds{1}_{[0,t]} (G_n(Y_{N(k)})) - \mathds{1}_{[0,t]} (G(Y_{N(k)})) \big| 
	  =: I_{1,s} + I_{2,t}.
\end{eqnarray*}
Since $G$ is continuous, we first obtain
\begin{eqnarray*}
  \lefteqn{\frac{1}{n} \, \sum_{k=1}^{n} 
						\big| \mathds{1}_{(-\infty,G_n^\leftarrow(s)]} (Y_{k}) - \mathds{1}_{(-\infty,G^\leftarrow(s)]} (Y_{k}) \big|}
	\\
	&   =  & \frac{1}{n} \, \sum_{k=1}^{n} \mathds{1}_{(-\infty,G_n^\leftarrow(s)]} (Y_{k})
					 + \frac{1}{n} \, \sum_{k=1}^{n} \mathds{1}_{(-\infty,G^\leftarrow(s)]} (Y_{k})
					 - \, 2 \, \frac{1}{n} \, \sum_{k=1}^{n} \mathds{1}_{(-\infty,G_n^\leftarrow(s)]} (Y_{k})\mathds{1}_{(-\infty,G^\leftarrow(s)]} (Y_{k})
	\\
	&   =  & | G_n \circ G_n^\leftarrow (s) - G_n \circ G^\leftarrow (s) |
	  \leq   | G_n \circ G_n^\leftarrow (s) - s | + | G \circ G^\leftarrow (s) - G_n \circ G^\leftarrow (s) |.
\end{eqnarray*}
By Glivenko-Cantelli, the second term converges to $0$ almost surely and,
for $(l-1)/n < s \leq l/n$, the first term reduces to  
$ | G_n \circ G_n^\leftarrow (s) - s |
	= | l/n - s | 
	\leq 1/n $.
Therefore,
\begin{eqnarray*}
  \lim_{n \to \infty} E(I_{1,s}) 
	&   =  & \lim_{n \to \infty} E \left( \frac{1}{n} \, \sum_{k=1}^{n} 
						\big| \mathds{1}_{[0,s]} (G_n(Y_k)) - \mathds{1}_{[0,s]} (G(Y_k)) \big| \right)
	\\*
	& \leq & \lim_{n \to \infty} \left(\frac{1}{n} + E \left( \frac{1}{n} \, \sum_{k=1}^{n} 
						\big| \mathds{1}_{(-\infty,G_n^\leftarrow(s)]} (Y_{k}) - \mathds{1}_{(-\infty,G^\leftarrow(s)]} (Y_{k}) \big| \right) \right)
	   =     0.
\end{eqnarray*}
Now, define $f(Y_{N(k)}) := \big| \mathds{1}_{[0,t]} (G_n(Y_{N(k)})) - \mathds{1}_{[0,t]} (G(Y_{N(k)})) \big|$.
Due to Lemma \ref{Lemma.Chatterjee.1} and the fact that ties only occur with probability $0$, we obtain 
$$
	\frac{1}{n} \, \sum_{k=1}^{n} E(f(Y_{N(k)}))
	  \leq \frac{1}{n} \, \sum_{k=1}^{n} \sum_{i=1}^{n} E\big( f(Y_{i}) \, \mathds{1}_{\{Y_i=Y_{N(k)}\}} \big)
	    =  \frac{1}{n} \, \sum_{i=1}^{n} E\left( f(Y_{i}) \, \sum_{k=1}^{n} \mathds{1}_{\{Y_i=Y_{N(k)}\}} \right)
	  \leq \frac{1}{n} \, \sum_{i=1}^{n} E\big( f(Y_{i}) \, c(d) \big)
	  \leq c(d) \, E\big( f(Y_{1}) \big),
$$
and hence 
$E(I_{2,t}) 
  = \frac{1}{n} \, \sum_{k=1}^{n} E(f(Y_{N(k)})) 
	\leq c(d) \, E\big( f(Y_{1}) \big) 
	= c(d) \, \frac{1}{n} \, \sum_{k=1}^{n} E(f(Y_{k}))
	= c(d) \, E(I_{1,t})$
which implies \linebreak $ \lim_{n \to \infty} E(I_{2,t}) = 0 $ and eventually
$ \lim_{n \to \infty} E \big( | C_n (s,t)  - C_n^\ast (s,t) | \big) = 0 $.
This proves 1.
\\
In order to prove key result 2. we adapt the ideas used in the proof of \cite[Lemma 11.8]{chatterjee2021}.
Let $\mathcal{F}$ be the $\sigma$-field generated by $\XXX_1, \dots,\XXX_n$ and the random variables used for breaking ties when selecting nearest neighbours.
Then, for $(s,t) \in \I^2$, we have
\begin{eqnarray*}
  E \big( \mathds{1}_{[0,s]} (G(Y_1)) \, \mathds{1}_{[0,t]} (G(Y_{N(1)})) \big| \mathcal{F} \big)
	& = & E \big( \mathds{1}_{[0,s]} (G(Y_1)) \big| \mathcal{F} \big) \; E \big( \mathds{1}_{[0,t]} (G(Y_{N(1)})) \big| \mathcal{F} \big)
	\\
	& = & P \big( G(Y_1) \leq s \big| \XXX_1 \big) \; P \big( G(Y_{N(1)}) \leq t \big| \XXX_{N(1)} \big).
	%\\
	%& = & P \big( F(Y) \leq s | \XXX_1 \big) P \big( F(Y) \leq t | \XXX_{N(1)} \big)
\end{eqnarray*}
Further, let $f_t(\XXX) := P \big( G(Y) \leq t | \XXX \big)$.
Then $f_t$ is measurable and, by \cite[Lemma 11.7]{chatterjee2021}, $f_t(\XXX_{N(1)}) - f_t(\XXX_1)$ tends to $0$ in probability,
and since $|f_s(\XXX_1)| = P \big( G(Y_1) \leq s \big| \XXX_1 \big) \leq 1$ it follows that (see, e.g., \cite{klenke2008})
\begin{eqnarray*}
  \lefteqn{\lim_{n \to \infty} E \big( \mathds{1}_{[0,s]} (G(Y_1)) \, \mathds{1}_{[0,t]} (G(Y_{N(1)})) \big)
	  =   \lim_{n \to \infty} E \left( E \big( \mathds{1}_{[0,s]} (G(Y_1)) \, \mathds{1}_{[0,t]} (G(Y_{N(1)})) \big| \mathcal{F} \big) \right)}
	\\
	& = & \lim_{n \to \infty} E \left( P \big( G(Y_1) \leq s | \XXX_1 \big) \,  P \big( G(Y_{N(1)}) \leq t \big| \XXX_{N(1)} \big) \right)
	  =   E \left( P \big( G(Y_1) \leq s \big| \XXX_1 \big) \; P \big( G(Y_1) \leq t \big| \XXX_1 \big) \right)
	  =   \psi(A) (s,t), 
\end{eqnarray*}
and hence $\lim_{n \to \infty} E (C_n^\ast(s,t)) = \psi(A) (s,t)$.
This proves 2. from which the assertion follows.
\qed

\bigskip
The proof of the following lemma is similar to that of \cite[Lemma 11.9]{chatterjee2021}:

\begin{lemma}{} \label{Consistency.L2} %{\rm [compare Azadkia and Chatterjee, Lemma 11.9]} \\
There are positive constants $M_1$ and $M_2$ depending only on dimension $d$ such that for any $n \in \mathbb{N}$ and any 
$\eta \in (0,\infty)$
$$ P (\{|C_n(s,t) - E(C_n(s,t))| \geq \eta\})
	\leq M_1 \exp(-M_2n\eta^2), 
	\qquad (s,t) \in \I^2.
$$
\end{lemma}{}
\noindent
{\bf Proof.} 
We prove the result using the bounded difference concentration inequality (see, e.g., \cite{mcdiarmid1989}).
For that purpose, we consider a realization $(\XXX_1, Y_1), \dots, (\XXX_n, Y_n)$ and i.i.d. uniformly ${\rm U}[0,1]$-distributed random variables $U_1, \dots, U_n$ where $U_i$ is used for breaking ties when selecting $\XXX_i$'s nearest neighbour.
We replace $(\XXX_i, Y_i, U_i)$ by some alternative value $(\XXX_i^\circ, Y_i^\circ, U_i^\circ)$ and find a bound on the maximum possible change in $C_n$.
After such a replacement, 
each $G_n(Y_j)$, $j \neq i$, can change by at most $1/n$.
\\
Now, 
for $(s,t) \in \I^2$, define
$ a_{k} 
	:= \mathds{1}_{[0,s]} (G_n(Y_k)) $ and 
$	b_{k} 
	:= \mathds{1}_{[0,t]} (G_n(Y_k)) $.
Then 
$ C_n (s,t) 
	=	\frac{1}{n} \sum_{k=1}^{n} a_k \, b_{N(k)} $.
When replacing $Y_i$ by $Y_i^\circ$ there exist at most two indices $l$ such that $a_l$ changes by $1$ 
and, by Lemma \ref{Lemma.Chatterjee.1}, there are at most $2 \cdot c(d)$ indices $l$ such that $b_{N(l)}$ changes by $1$.
When replacing $\XXX_i$ by $\XXX_i^\circ$, by Lemma \ref{Lemma.Chatterjee.1}, there are at most $2 \cdot c(d)+1$ indices $l$ such that $b_{N(l)}$ changes by $1$.
And when replacing $U_i$ by $U_i^\circ$ there exist at most one index $l$ such that $b_{N(l)}$ changes by $1$.
In sum, there are at most $4 \cdot c(d) +4$ indices $l$ such that $a_l \, b_{N(l)}$ changes by $1$,
and hence $C_n (s,t)$ changes by at most $(4 \cdot c(d) +4)/n$.
Applying the bounded difference concentration inequality and replacing $(4 \cdot c(d) +4)^{-2}$ by $c(d)$, 
for every $\eta \in (0,\infty)$, we obtain
$ P (\{|C_n(s,t) - E(C_n(s,t))| \geq \eta\})
	\leq 2 \, \exp(- c(d)n\eta^2) $.
As the above inequality holds for all $(s,t) \in \I^2$, this completes the proof.
\qedsymbol

\bigskip
Combining Lemma \ref{thm.convergence.1} and Lemma \ref{Consistency.L2} proves consistency of $C_n$:

\begin{lemma}{} \label{Consistency.Thm}
Let $(\XXX_1,Y_1), (\XXX_2,Y_2),\ldots$ be a random sample from $(\XXX,Y)$ with continuous marginals and connecting copula $A \in \mathcal{C}^{d+1}$.
Then, for all $(s,t) \in \I^2$, we have
$ \lim_{n \to \infty} C_n (s,t) = \psi(A) (s,t) $
almost surely.
\end{lemma}{}

Due to Lemma \ref{Lemma.Chatterjee.1}, we have
\begin{equation} \label{Est.Asymtotic.Eq}
  | C_n(s,t) - D_n(s,t) | \leq \frac{c(d)+1}{n},
\end{equation}
implying that $C_n$ and $D_n$ are asymptotically equivalent. This proves consistency of $D_n$ and hence Theorem \ref{Consistency.Thm2}.
\end{document}